\documentclass[11pt,a4paper]{article}

\usepackage{graphicx}
\usepackage{amsfonts,amsmath, amssymb}

\numberwithin{equation}{section}

\usepackage{graphics}
\usepackage{epsfig}
\allowdisplaybreaks

\newtheorem{claim}{\bf \t}[part]

\newtheorem{theorem}{Theorem}[section]

\newtheorem{lemma}[theorem]{Lemma}
\newtheorem{proposition}[theorem]{Proposition}
\newtheorem{remark}[theorem]{Remark}

\def\t{\theta}

\begin{document}

\title{Planar Traveling Waves For
 Nonlocal Dispersion Equation With Monostable Nonlinearity}

\author{ Rui Huang$^a$, \ Ming Mei$^{b,c}$
  \ and \
 Yong Wang$^{d}$
\\
\ \\
   {\small \it $^a$School of Mathematical Sciences, South China Normal University}\\
   {\small \it Guangzhou, Guangdong, 510631, China}\\
   {\small \tt huang@scnu.edu.cn}\\
   {\small and }\\
   {\small \it $^b$Department of Mathematics, Champlain College Saint-Lambert} \\
   {\small\it Quebec,  J4P 3P2, Canada } \\
  {\small \it $^c$Department of Mathematics and Statistics, McGill
University} \\
{\small \it  Montreal, Quebec,  \ H3A 2K6, Canada} \\
{\small \tt ming.mei@mcgill.ca}\\
 {\small and} \\
   {\small \it $^d$Institute of Applied Mathematics,  Academy of Mathematics and Systems
   Science}\\
{\small \it Chinese Academy of Sciences, Beijing, 100190, China} \\
{\small \tt yongwang@amss.ac.cn}\\
}
\date{ }
\maketitle

\begin{abstract}
In this paper, we study a class of  nonlocal dispersion
equation with monostable nonlinearity in $n$-dimensional space
\[
\begin{cases}
 u_t - J\ast u +u+d(u(t,x))=
\int_{\mathbb{R}^n} f_\beta
(y) b(u(t-\tau,x-y)) dy, \\
u(s,x)=u_0(s,x), \ \ s\in[-\tau,0], \ x\in \mathbb{R}^n,
\end{cases}
\]
where the nonlinear functions $d(u)$ and $b(u)$ possess the monostable characters like Fisher-KPP type, $f_\beta(x)$ is the heat kernel, and the kernel $J(x)$ satisfies
 ${\hat J}(\xi)=1-\mathcal{K}|\xi|^\alpha+o(|\xi|^\alpha)$ for $0<\alpha\le 2$. After establishing the existence
  for both the planar traveling waves $\phi(x\cdot{\bf e}+ct)$ for $c\ge c_*$ ($c_*$ is the critical wave speed) and the solution $u(t,x)$ for the Cauchy problem, as well as the comparison principles, we prove that,
 all noncritical planar wavefronts $\phi(x\cdot{\bf e}+ct)$ are globally stable with
the exponential convergence rate $t^{-n/\alpha}e^{-\mu_\tau}$ for $\mu_\tau>0$, and the critical
wavefronts $\phi(x\cdot{\bf e}+c_*t)$ are globally stable in the algebraic form $t^{-n/\alpha}$. The adopted approach is  Fourier transform
and the weighted energy
method  with a suitably selected weight function. These rates are optimal and the stability results significantly develop the existing studies for nonlocal dispersion equations.

\


\

AMS: 35K57, 34K20, 92D25
\end{abstract}

\newpage

\tableofcontents


\baselineskip=16pt
\section{Introduction}\label{int}

For the gradient flow to an order parameter describing the state of a solid material, for example, a perfect
crystal with two different orientations, it is usually described by a convolution model of phase transition in the form \cite{Bates-Fife-Ren-Wang,Chasseigne-Chaves-Rossi, Cortazar-Elgueta-Rossi-Wolanski, Ignat-Rossi-1,Ignat-Rossi-2}
\begin{equation}
u_t = J\ast u-u +F(u), \ \ (x,t)\in \mathbb{R}^n\times \mathbb{R}_+,
\label{01}
\end{equation}
where $x=(x_1,x_2,\cdots,x_n)\in \mathbb{R}^n$,   $J(x)$  is a non-negative and radial kernel with unit integral, and
\begin{equation}
(J\ast u) (t,x)=\int_{\mathbb{R}^n} J(x-y) u(t,y) dy.
\label{02}
\end{equation}
 As showed in \cite{Chasseigne-Chaves-Rossi,Coville-Dupaigne-1},
when the kernel $J(x)$ has a second momentum, for example, $J$ is compact-supported or Gaussian-like kernel $J\sim e^{-x^2}$,  its Fourier transform looks like
\[
{\hat J}(\xi)=1-\mathcal{K}|\xi|^2 + o(|\xi|^2), \ \ \mathcal{K}>0,
\]
then the effect of the nonlocal dispersion $J\ast u-u$ is
almost the same to the linear diffusion $\mathcal{K}\Delta u$:
\[
 J\ast u-u \ \ \approx  \ \ \mathcal{K}\Delta u,
 \]
which informs us to expect that the behaviors of the solutions to the nonlocal dispersion equation and the linear diffusion equation are almost
identical \cite{Chasseigne-Chaves-Rossi, Cortazar-Elgueta-Rossi-Wolanski, Ignat-Rossi-1, Ignat-Rossi-2}
\[
u_t= J\ast u-u \ \ \Leftrightarrow \ \  u_t=\mathcal{K} \Delta u.
\]
Notice that, comparing with the heat equations, the solutions for the nonlocal dispersion equations usually loss the spatial regularity,
but have much better regularity in time, see Remark \ref{remark} below for details.

In general, $J(x)$ may not  have a second momentum, let us say,
\[
{\hat J}(\xi)= 1-\mathcal{K}|\xi|^\alpha +o(|\xi|^\alpha) \ \
\mbox{as}\  \xi\rightarrow0 \ \mbox{ for } \ \alpha \in (0, 2).
\]
One example is the Cauchy law by taking $J(x)=\frac{1}{1+|x|^2}$ which implies its Fourier transform mentioned above with $\alpha=1$. In this case, the behavior of the solutions to the nonlocal dispersion equation is almost identical to
the fractional diffusion equation \cite{Chasseigne-Chaves-Rossi,Cortazar-Elgueta-Rossi-Wolanski, Ignat-Rossi-1, Ignat-Rossi-2}
\[
 u_t= J\ast u-u \ \ \Leftrightarrow \ \  u_t=\mathcal{K} \Delta^{\alpha/2} u.
\]

Equation \eqref{01} represents  also the dynamical population model of single species in ecology \cite{Fife}, where $u(t,x)$ is the density of population at location $x$ and time $t$, and $J(x-y)$ is thought of as the probability distribution of jumping from location $y$ to location $x$, and $J*u=\int_{\mathbb{R}^n}J(x-y) u(t,y) dy$ is the rate at which individuals are arriving to position $x$ from all other places, while $-u(x,t)=-\int_{\mathbb{R}^n}J(x-y) u(t,x) dy$ stands the rate at which they are leaving the location $x$ to travel to all other places. In this case, under the consideration of the effects from  birth rate and death rate, the equation \eqref{01} is usually written as follows
\begin{equation}
u_t = J\ast u -u +b(u(t-\tau,x))-d(u(t,x)), \ \ (x,t)\in \mathbb{R}^n\times \mathbb{R}_+,
\label{03}
\end{equation}
where $b(u(t-\tau, x)$ is the birth rate function,  $d(u(t,x))$ is the death rate function, and $\tau>0$ is the mature age of the single species, which is usually called the {\it time-delay}. Furthermore, if we consider the distribution of all matured population, the effect of birth rate is then involved in whole space $\mathbb{R}^n$ \cite{Gourley-Wu,Mei-Ou-Zhao,So-Wu-Zou}, and the equation is expressed as
\begin{equation}
\frac{\partial u}{\partial t} - J\ast u +u+d(u(t,x))=
\int_{\mathbb{R}^n} f_\beta
(y) b(u(t-\tau,x-y)) dy,
\label{04}
\end{equation}
where
$f_\beta(y)$, with $\beta>0$,  is the heat kernel in the form of
\begin{equation}
f_\beta (y)=\frac{1}{(4\pi \beta)^\frac n2}
e^{\frac{-|y|^2}{4\beta}} \ \ \mbox{ with } \ \ \int_{\mathbb{R}^n}
f_\beta (y) dy =1. \label{1.2}
\end{equation}
Notice that, by using the property of heat kernel
\[
\lim_{\beta\to 0^+} \int_{\mathbb{R}^n} f_\beta
(y) b(u(t-\tau,x-y)) dy= b(u(t-\tau,x)),
\]
we then derive the equation \eqref{03}  as a limit of the equation \eqref{04}  by taking $\beta\to 0^+$, and further derive the regular nonlocal dispersion equation \eqref{01} from the equation \eqref{03} by taking the time-delay $\tau=0$ and $F(u)=b(u)-d(u)$.  In particular, if we set $d(u)=u^2$ and $b(u)=u$,  then, from \eqref{01} we get the classical Fisher-KPP equation with nonlocal dispersion
\begin{equation}
u_t=J\ast u-u + u(1-u). \label{05}
\end{equation}
So, the equations \eqref{01} and \eqref{03} and \eqref{05}  all are the special cases
of the equation \eqref{04}.

In this paper, we will concentrate ourselves to  the Cauchy problem for the more generalized equation \eqref{04} with non-locality of birth rate
\begin{equation}
\begin{cases}
\displaystyle \frac{\partial u}{\partial t} - J\ast u +u+d(u(t,x))=
\int_{\mathbb{R}^n} f_\beta
(y) b(u(t-\tau,x-y)) dy, \\
u(s,x)=u_0(s,x), \ \ s\in[-\tau,0], \ x\in \mathbb{R}^n.
\end{cases}
\label{1.1}
\end{equation}
When $\tau=0$ (no time-delay), then the above equation is reduced to
\begin{equation}
\begin{cases}
\displaystyle \frac{\partial u}{\partial t} - J\ast u +u+d(u)=
\int_{\mathbb{R}^n} f_\beta
(y) b(u(t,x-y)) dy, \\
u(0,x)=u_0(x), \  \ x\in \mathbb{R}^n.
\end{cases}
\label{1.1-2}
\end{equation}
We will also discuss how the time-delay $\tau$ effects the property of the solutions.

For the equation \eqref{01} in 1D case, when $F(u)$ is bistable, namely, two constant equilibria $u_-$ and $u_+$ both are the stable nodes
(the typical example   is the Huxley equation with  $F(u)=u(u-a)(1-u)$ for $0<a<1$),
Bates {\it et al} \cite{Bates-Fife-Ren-Wang} and Chen \cite{XChen} proved that the traveling waves are globally stable as $t\to+\infty$. In this paper, we consider another important  type of equations
with monostable nonlinearity. The typical example in this case is  Fisher-KPP equation with  $F(u)=u(1-u)$. Hence, throughout this paper, we assume that
the death rate $d(u)$ and birth rate $b(u)$ capture the following characters of monostable nonlinearity:
\begin{enumerate}
\item[(H$_1$)]  There exist $u_-=0$ and $u_+>0$ such that
$d(0)=b(0)=0$,  $d(u_+)= b(u_+)$, and $d(u),b(u)\in C^2[0,u_+]$;

\item[(H$_2$)]  $ b'(0)>d'(0)\ge 0$ and $0\le
b'(u_+)<d'(u_+)$;

\item[(H$_3$)] For $0\le u\le u_+$, $d'(u)\ge 0$, $b'(u)\geq 0$,
$d''(u)\ge 0$, $b''(u)\le 0$.
\end{enumerate}
These characters are summarized from the classical Fisher-KPP equation, see also the monostable reaction-diffusion equations in ecology, for example, the Nicholson's blowflies equation  \cite{MLLS1,MLLS2,Mei-Ou-Zhao,So-Wu-Zou} with
\[
d(u)=\delta u  \ \mbox{ and } \  b(u)= p ue^{-a u}, \  p>0, \delta>0, a>0
\]
and $u_-=0$ and $u_+=\frac{1}{a}\ln \frac{p}{\delta}>0$ under the consideration of $1<\frac{p}{\delta}\le e$; and the age-structured population model
\cite{Gourley2,Gourley-Wu,Mei-Ou-Zhao,Mei-Wong,Pan-Li-Lin} with
\[
d(u)=\delta u^2 \ \mbox{ and } \  b(u)=pe^{-\gamma\tau}u, \ \ \delta>0, \ p>0, \ \gamma>0,
\]
and $u_-=0$ and $u_+=\frac{p}{\delta}e^{-\gamma\tau}$.

Clearly, under the hypothesis (H$_1$)-(H$_3$), both $u_-=0$ and $u_+>0$ are
constant equilibria of the equation \eqref{1.1}, and
 $u_-=0$ is unstable and $u_+$ is stable for the
spatially homogeneous equation associated with \eqref{1.1}, this is why we call the equation \eqref{1.1}, including \eqref{01} and
\eqref{03} and \eqref{1.1-2}, as monostable.

On the other hand,
we also assume the kernel $J(x)$ satisfying:
\begin{enumerate}
\item[(J$_1$)]  $\displaystyle J(x)=\prod^n_{i=1}J_i(x_i)$, where $J_i(x_i)$  is smooth, and $J_i(x_i)=J_i(|x_i|)\ge 0$ and $\displaystyle \int_\mathbb{R} J_i(x_i) dx_i=1$ for $i=1,2\cdots,n$, and $\int_{\mathbb{R}}|y_1|J_1(y_1)e^{-\lambda_* y_1} dy_1<\infty$ for $\lambda_*>0$ defined in \eqref{2.4''} and \eqref{2.4'};

\item[(J$_2$)]  Fourier transform of $J(x)$ satisfies ${\hat J}(\xi)=1-\mathcal{K}|\xi|^\alpha + o(|\xi|^\alpha)$ \mbox{as} $\xi\rightarrow0$\  with $\alpha \in (0,2]$ and $\mathcal{K}>0$.
\end{enumerate}

A {\it planar traveling wavefront} to the equation \eqref{1.1} for $\tau\ge 0$  is a special solution in the form of
$u(t,x)=\phi(x\cdot {\bf e}+ct)$ with $\phi(\pm\infty)=u_\pm$, where
$c$ is the wave speed, ${\bf e}$ is a unit vector of the basis of
$\mathbb{R}^n$. Without loss of generality, we can always assume
${\bf e}={\bf e}_1=(1,0,\cdots,0)$ by rotating the coordinates.
Thus,  the planar traveling wavefront
$\phi(x\cdot {\bf e}_1+ct)=\phi(x_1+ct)$ satisfies, for
$\tau\ge 0$,
\begin{equation}
\begin{cases}
c\phi'-J\ast\phi+\phi+d(\phi)=\displaystyle\int_{\mathbb{R}^n}f_\beta(y)
b(\phi(\xi_1-y_1-c\tau))dy, \\
\phi(\pm\infty)=u_\pm,
\end{cases}
\label{2.2}
\end{equation}
where $'=\frac{d}{d\xi_1}$ and  $\xi_1=x_1+ct$.  Let
\begin{equation}
f_{i\beta}(y_i):=\frac{1}{(4\pi\beta)^{1/2}}e^{-\frac{y_i^2}{4\beta}}.
\label{2.2-2}
\end{equation}
Then
\begin{equation}
f_\beta(y):=\prod^n_{i=1} f_{i\beta}(y_i), \ \mbox{ and }
\int_{\mathbb{R}}f_{i\beta}(y_i)dy_i=1, \ \ i=1,2,\cdots,n,
\label{2.2-3}
\end{equation}
and \eqref{2.2} is reduced to, for $\tau\ge 0$,
\begin{equation}
\begin{cases}
c\phi'-J_1\ast\phi+\phi+d(\phi)=\displaystyle\int_{\mathbb{R}}f_{1\beta}(y_1)
b(\phi(\xi_1-y_1-c\tau))dy_1, \\
\phi(\pm\infty)=u_\pm.
\end{cases}
\label{2.2-4}
\end{equation}
 The main
purpose of this paper is to study the global asymptotic stability of
planar traveling wavefronts of the equations \eqref{1.1} and
\eqref{1.1-2} with or without time-delay, respectively, in particular, in the
case of the {\it critical wave} $\phi(x_1+c_*t)$. Here the number
$c_*$ is called the {\it critical speed} (or the {\it minimum
speed}) in the sense that a traveling wave $\phi(x_1+ct)$ exists if
$c\ge c_*$, while no traveling wave $\phi(x_1+ct)$ exists if
$c<c_*$.

The nonlocal dispersion equation \eqref{01} has been extensively studied recently. Chasseigne {\it et al} \cite{Chasseigne-Chaves-Rossi}
and Cortazar {\it et al} \cite{Cortazar-Elgueta-Rossi-Wolanski} showed that the linear nonlocal dispersion equation \eqref{01} (with $F=0$)
is almost equivalent to the linear diffusion equation, and the asymptotic behavior of the solutions to the linear  equation of nonlocal dispersion is exactly
the same to the corresponding linear diffusion equation. Ignat and Rossi \cite{Ignat-Rossi-1,Ignat-Rossi-2} further obtained the asymptotic behavior
of the solutions to the nonlinear equation \eqref{01}. Garc$\acute{\mbox{i}}$a-Meli$\acute{\mbox{a}}$n and Quir$\acute{\mbox{o}}$s \cite{Garcia-Quiros} investigated the blow up phenomenon of the solution to the equation \eqref{01} with $F(u)=u^p$, and gave the Fujita critical exponent. Regarding the structure of special solutions to \eqref{01} like traveling wave solutions, early in 1997 Bates {\it et al} \cite{Bates-Fife-Ren-Wang} and Chen \cite{XChen} established the existence of the traveling waves for \eqref{01} with bistable nonlinearity, and proved their global stability. For \eqref{01} with monostable nonlinearity,  recently
 Coville and his collaborators \cite{Coville,Coville-Davila-Martinez,Coville-Dupaigne-1,Coville-Dupaigne-2} studied the existence and uniqueness (up to a shift) of traveling waves.
See also the existence/nonexistence of traveling waves by Yagisita \cite{Y} and the existence of almost periodic traveling waves by Chen \cite{Chen}. However, the stability of traveling waves for the nonlocal equation \eqref{1.1} (including \eqref{01} and \eqref{03}) with monostable nonlinearity is almost not related, except a special case for the fast waves with large wave speed to the 1D age-structured population model by Pan {\it et al} \cite{Pan-Li-Lin}. As we know,  such a problem is also
very significant but challenging, because the  equations of Fisher-KPP type possess an unstable node, different from the bistable case, this unstable node usually causes a serious difficulty in the stability proof, particularly, for the critical traveling waves.  The main interest in this paper is to investigate the stability of traveling waves to \eqref{1.1} with $\tau>0$ and \eqref{1.1-2} with $\tau=0$. An  easy to follow method will be
introduced for the stability proof to the nonlocal dispersion equations.

In this paper, we will first investigate  the linearized equation of \eqref{1.1}, and derive the optimal decay rates of the solution to the linearized equation by means of
Fourier transform. This is a crucial step for get the optimal convergence for the nonlocal stability of traveling waves. Then, we will technically establish the global existence and comparison principles of the solution to the $n$-D
nonlinear equation with nonlocal dispersion \eqref{1.1}. Inspired by \cite{Meot} for the classical Fisher-KPP equations and the further developments by \cite{Mei-Ou-Zhao,Mei-Ou-Zhao2},  by ingeniously selecting a weight function which is dependent on the critical wave speed $c_*$,
 and using the weighted energy method and the Green function method with the comparison principles together, we will further prove that, all noncritical planar traveling waves $\phi(x\cdot {\bf e} +ct)$ are exponentially stable in the form of
$t^{-\frac{n}{\alpha}}e^{-\mu_\tau}$ for some constant $\mu_\tau=\mu(\tau)$ such that $0<\mu_\tau\le \mu_0$ for $\tau\ge 0$; and all critical planar traveling waves $\phi(x\cdot {\bf e} +c_*t)$ are algebraically stable in the form of $t^{-\frac{n}{\alpha}}$. These convergence rates are optimal and the stability results significantly develop the existing studies on the nonlocal dispersion equations.
We will also show that the time-delay $\tau$ will slow down  the convergence of the the solution $u(t,x)$ to the noncritical planar traveling waves
$\phi(x\cdot{\bf e}+ct)$ with $c>c_*$, and cause the higher requirement for the initial perturbation around the wavefronts.

For the stability of traveling waves to other modeling equations, we refer to the classical and significant contributions in
\cite{AW2,Bramson,Chen-Guo-Wu,FM,Ga,HR,Huang,Ka,K,KPP,Lau,MR,MNT,MLLS1,MLLS2,Mei-Ou-Zhao,Mei-Wong,Meot,S,SmZ,U,VVV,Xin,X} for reaction-diffusion equations and \cite{FS,Goodman,Kawashima-Matsumura-1,Kawashima-Matsumura-2,Matsumura-Mei,Matsumura-Nishihara-1,Matsumura-Nishihara-2,SX,ZS} for fluid dynamical systems, and the references therein.

The paper is organized as follows. In section 2, we will state the existence of the traveling waves, and their stability. In section 3, we will give the solution formulas to the linearized dispersion equations of \eqref{1.1} and \eqref{1.1-2}, and derive the optimal decay rates by Fourier transform with energy method together. In section 4, we will prove the global existence of the solution to \eqref{1.1} and establish the comparison principle.
In section 5, based on the results obtained in sections 3 and 4, by using the weighted energy method, we will further prove the stability of planar traveling waves including the critical and noncritical waves. Finally, in section 6, we will give some particular applications of our stability theory to the classical Fisher-KPP equation with nonlocal dispersion and the Nicholson's blowflies model,  and make a concluding remark to a more general case.

 Before ending this section, we make some notations. Throughout this paper, $C>0$ denotes a generic constant, while
$C_i>0$ and $c_i>0$ ($i=0,1,2,\cdots$)
 represent specific constants. $j=(j_1,j_2,\cdots, j_n)$ denotes a multi-index with non-negative integers $j_i\ge 0$ ($i=1,\cdots,n$), and
 $|j|=j_1+j_2+\cdots +j_n$.
 The derivatives for multi-dimensional function are denoted as
 \[
 \partial^j_x f(x):= \partial^{j_1}_{x_1}\cdots \partial^{j_n}_{x_n}f(x).
 \]
For a $n$-D function $f(x)$, its Fourier transform  is defined as
\[
\mathcal{F}[f](\eta)=\hat{f}(\eta):=\int_{\mathbb{R}^n}
e^{-\mbox{i}x\cdot \eta} f(x) dx, \ \ \ \ \ \mbox{i}:=\sqrt{-1},
\]
and the inverse Fourier transform is given by
\[
\mathcal{F}^{-1}[\hat{f}](x):=\frac{1}{(2\pi)^n} \int_{\mathbb{R}^n}
e^{\mbox{i}x\cdot \eta} \hat{f}(\eta) d\eta.
\]
 Let $I$ be an interval,
typically $I = \mathbb{R}^n$. $L^p(I)$ ($p\ge 1$) is the Lebesque
space of the integrable
  functions defined on  $I$, $W^{k,p}(I)$ ($k\ge 0, p\ge 1$)  is
the Sobolev space of the $L^p$-functions $f(x)$ defined on the
interval $I$ whose derivatives $\partial^j_x f$ with $|j|=k$ also
  belong to $L^p(I)$, and in particular, we denote
  $W^{k,2}(I)$ as $H^k(I)$. Further,
$L^p_w(I)$ denotes the  weighted $L^p$-space for a weight function
$w(x)>0$ with the norm defined as
\[
\|f\|_{L^p_w}=\Big(\int_{I} w(x) \left | f(x) \right | ^p
dx\Big)^{1/p},
\]
$W^{k,p}_w(I)$ is the weighted Sobolev space with the norm given by
\[
\|f\|_{W^{k,p}_w}=\Big( \sum_{|j|=0}^k \int_{I} w(x) \left |
\partial^j_x f(x)\right |^p dx\Big)^{1/p},
\]
and $H^k_w(I)$ is defined with the norm
\[
\|f\|_{H^k_w}=\Big( \sum_{|j|=0}^k \int_{I} w(x) \left |
\partial^j_x f(x)\right | ^2 dx\Big)^{1/2}.
\]
 Let $T > 0$ be a number and  ${\cal B}$ be
a Banach space. We denote by $C^0([0,T],{\cal B})$ the space of the
$\cal B$-valued continuous functions on $[0,T]$, $L^2([0,T], {\cal
B})$ as the space of the ${\cal B}$-valued $L^2$-functions on
$[0,T]$. The corresponding spaces of the ${\cal B}$-valued functions
on $[0,\infty)$ are defined similarly.

\section{Traveling Waves and Their Stabilities}

As we mentioned before, the existence and uniqueness (up to a shift)
of traveling waves for the equation \eqref{01} were proved in
\cite{Coville,Coville-Davila-Martinez,Coville-Dupaigne-1,Coville-Dupaigne-2},
particular, in a recent work by Yagisita \cite{Y} for the existence
and nonexistence of traveling waves, when the nonlinearity $F(u)$ is
monostable. Without any difficulty, these results can be extended to
the nonlocal equation \eqref{1.1} with time-delay with the help of
comparison principle established in Section 4, when $d(u)$ and
$b(u)$ satisfy the monostable features (H$_1$)-(H$_3$). We state
these results as follows  without detailed proof.

\begin{theorem}
 Under the conditions (H$_1$)-(H$_3$) and (J$_1$)-(J$_2$), for the time-delay $\tau\ge 0$, there exist a minimum
 wave speed (also called the critical wave speed)
 $c_*>0$ such that
 \begin{enumerate}
 \item[$\bullet$] when $c\ge c_*$, there exits a monotone
traveling wavefront $\phi(x_1+ct)$ of \eqref{2.2} connecting $u_\pm$
exists;
\item[$\bullet$] when $c<c_*$, no traveling wave $\phi(x_1+ct)$
exists.
\end{enumerate}
 Here $(c_*,\lambda_*)$ with $c_*>0$ and $\lambda_*>0$ is given by
 \begin{equation}
H_{c_*}(\lambda_*)=G_{c_*}(\lambda_*), \ \ \ \
H'_{c_*}(\lambda_*)=G'_{c_*}(\lambda_*), \label{2.3}
\end{equation}
where
\begin{equation}
H_{c}(\lambda)= b'(0) e^{\beta\lambda^2-\lambda c \tau}, \
G_{c}(\lambda)=c\lambda-E_c(\lambda) +d'(0),\
E_c(\lambda)=\int_{\mathbb{R}}J_1(y_1)e^{-\lambda y_1}dy_1-1,
\label{2.4}
\end{equation}
namely, $(c_*,\lambda_*)$ is the tangent point of $H_{c}(\lambda)$ and
$G_{c}(\lambda)$ specified as
\begin{eqnarray}
 b'(0)
e^{\beta\lambda_*^2-\lambda_* c_* \tau}&=&c_*\lambda_* -\int_{\mathbb{R}}J_1(y_1)e^{-\lambda_* y_1} dy_1 +1 + d'(0),\label{2.4''} \\
 b'(0)
(2\beta\lambda_*-c_*\tau)e^{\beta\lambda_*^2-\lambda_*c_*\tau}&=&c_* + \int_{\mathbb{R}}y_1 J_1(y_1)e^{-\lambda_* y_1} dy_1.
\label{2.4'}
\end{eqnarray}
Furthermore, it can be verified:
\begin{enumerate}
\item[$\bullet$] In the case of $c>c_*$, there exist two numbers depending on $c$:
$\lambda_{1}=\lambda_1(c)>0$ and $\lambda_{2}=\lambda_2(c)>0$ as the
solutions to the equation $H_c(\lambda_{i})=G_c(\lambda_{i})$, i.e.,
\begin{equation}
 b'(0)
e^{\beta\lambda_{i}^2-\lambda_{i} c \tau}= c\lambda_{i}
-\int_{\mathbb{R}}J_1(y_1)e^{-\lambda_i  y_1} dy_1 + d'(0), \ \ \ i=1,2, \label{2.4'-new}
\end{equation}
such that
\begin{equation}
H_c(\lambda)<G_c(\lambda) \ \ \ \mbox{ for } \
\lambda_{1}<\lambda<\lambda_{2}, \label{2.5}
\end{equation}
and particularly,
\begin{equation}
H_c(\lambda_*)<G_c(\lambda_*) \ \mbox{ with } \
\lambda_{1}<\lambda_*<\lambda_{2}. \label{2.4'-newnew}
\end{equation}

 \item[$\bullet$] In the case of  $c=c_*$, it holds
\begin{equation}
H_{c_*}(\lambda_*)=G_{c_*}(\lambda_*) \ \mbox{ with } \
\lambda_{1}=\lambda_*=\lambda_{2}. \label{2.5-new}
\end{equation}

\item[$\bullet$] When $\xi_1=x_1+ct \to \pm\infty$,  for all $c\ge c_*$, the
traveling wavefronts $\phi(x_1+ct)$ converge to $u_\pm$
exponentially as follows
\begin{equation}\label{newnew}
|\phi(\xi_1)-u_\pm|=O(1) e^{-\lambda^\pm |\xi_1|}.
\end{equation}
Here $\lambda^{-}=\lambda_1(c)>0$ is given in \eqref{2.4'-new}, and
$\lambda^+=\lambda^+(c)>0$ is the unique root determined by the following
equation
\begin{equation}\label{newnewnew}
-c\lambda^+-\int_{\mathbb{R}}J_1(y_1)e^{-\lambda^+ y_1} dy_1+d'(u_+)=b'(u_+) e^{\beta
(\lambda^+)^2-\lambda^+ c\tau}.
\end{equation}
\end{enumerate}
 \label{TW}
 \end{theorem}

\begin{figure}
\centering
(a)\includegraphics[width=0.28\textwidth]{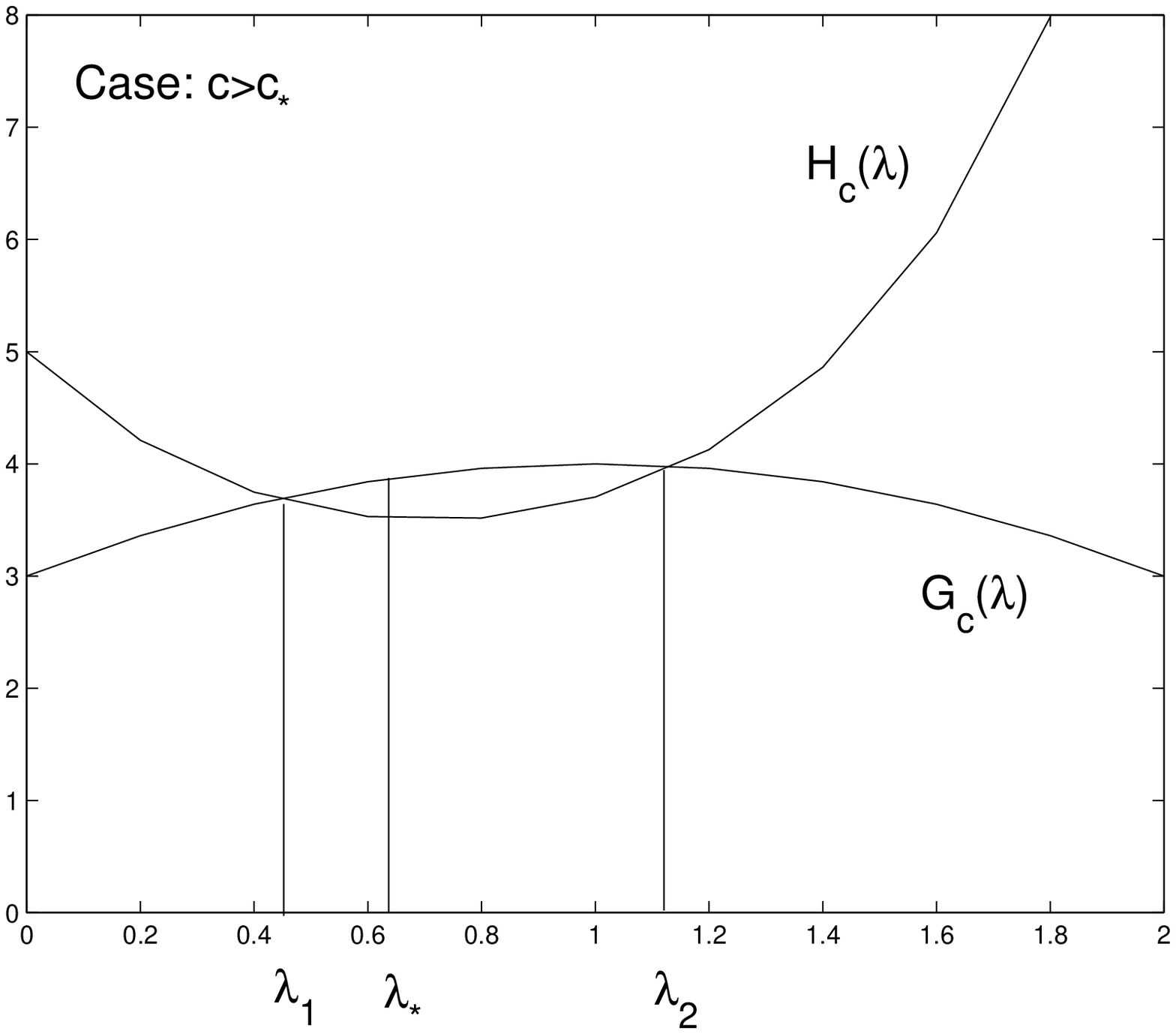}
  (b)\includegraphics[width=0.28\textwidth]{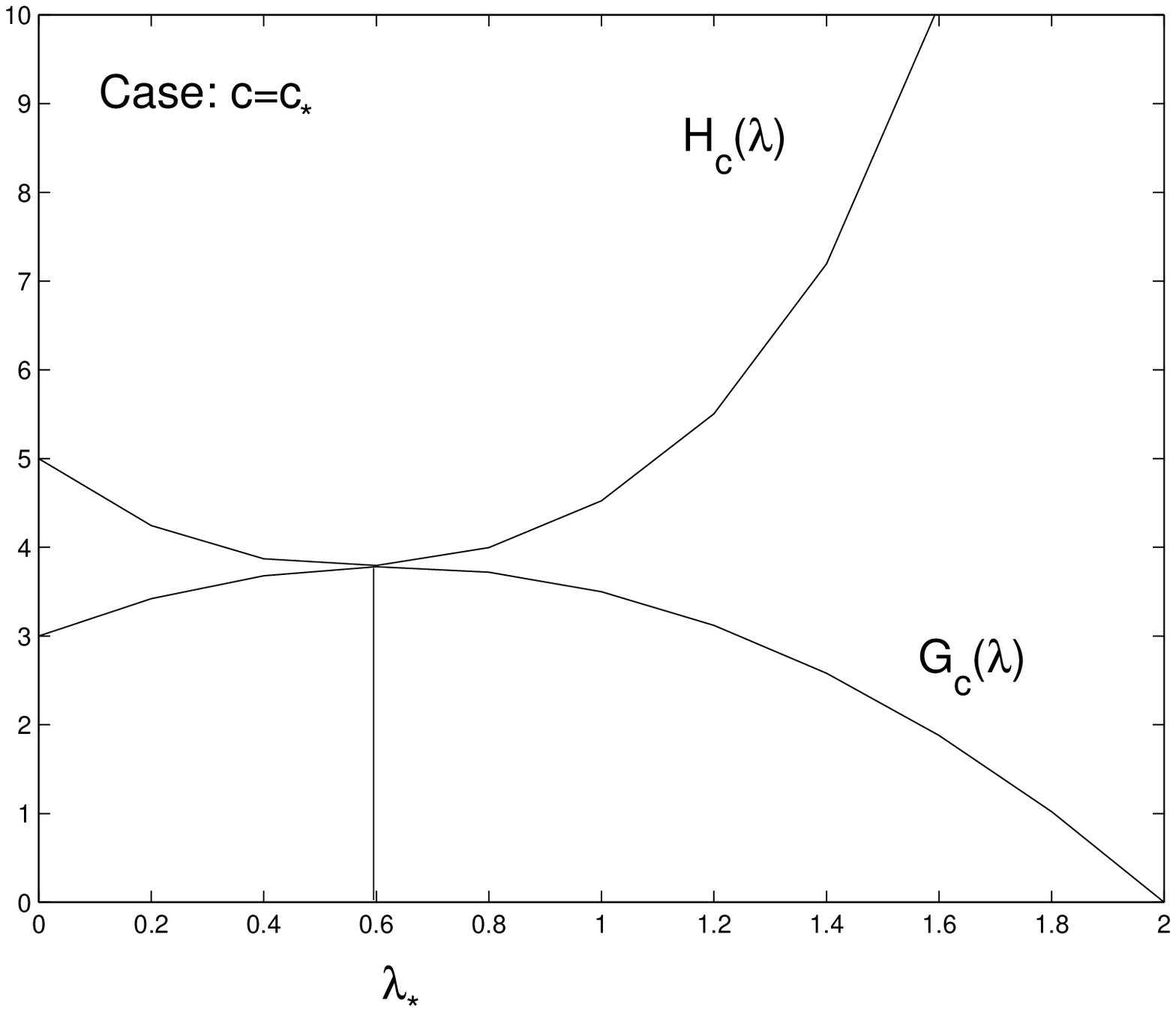}
  (c)\includegraphics[width=0.28\textwidth]{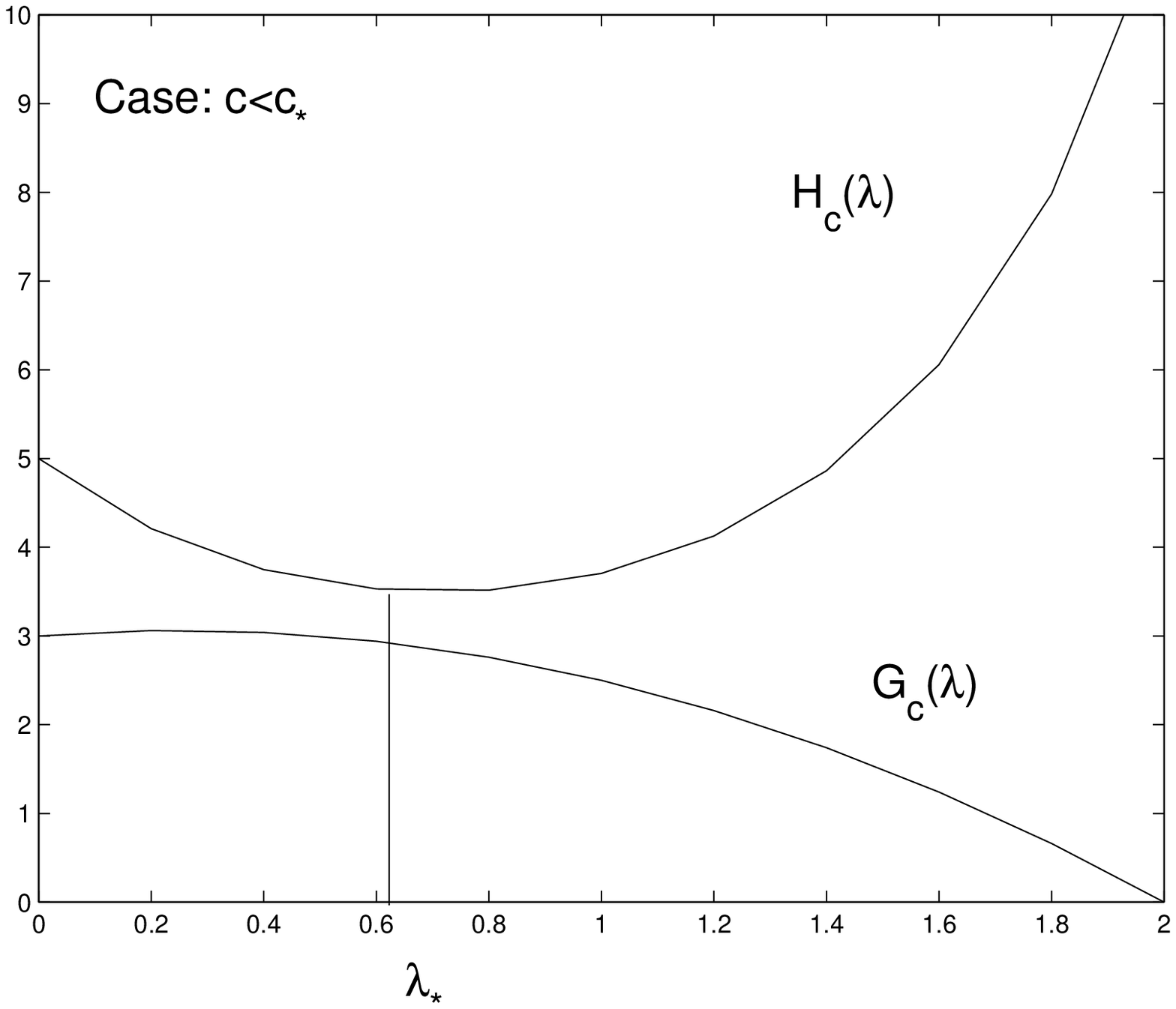}
  \caption{(a): the case of $c>c_*$; (b): the case of $c=c_*$; and (c): the case of $c<c_*$.}
  \label{fig}
\end{figure}

For easily understanding all cases mentioned in the above, we show them in Figure \ref{fig}.

Before stating our main stability theorems, let us technically choose a weight function:
\begin{equation}
w(x_1)=\begin{cases} e^{-\lambda_* (x_1-x_*)}, & \ \mbox{ for }
x_1\le
x_*, \\
1, & \ \mbox{ for } x_1>x_*,
\end{cases}
 \label{2.8-new}
\end{equation}
where $\lambda_*=\lambda_*(c_*)>0$ is given in \eqref{2.4''} and \eqref{2.4'}, and $x_*>0$ is a sufficiently large number such that,
\begin{equation}
0<d'(\phi(x_*))-\int_{\mathbb{R}^n} f_{\beta}(y)
b'(\phi(x_*-y_1-c\tau)) dy< d'(u_+)-b'(u_+).
\label{nov4-new}
\end{equation}
The selection of $x_*$  in \eqref{nov4-new} is valid, because of $d'(u_+)-b'(u_+)>0$ (see(H$_2$)). In fact, we have
\begin{eqnarray*}
\lim_{\xi_1\to \infty} d'(\phi(\xi_1))&=&d'(u_+)\\
&>&b'(u_+)\\
&=& \int_{\mathbb{R}^n} f_\beta(y)\Big[\lim_{\xi_1\to\infty} b'(\phi(\xi_1-y_1-c\tau))\Big] dy\\
&=&\lim_{\xi_1\to\infty}\int_{\mathbb{R}^n} f_\beta(y)
b'(\phi(\xi_1-y_1-c\tau))dy,
\end{eqnarray*}
which implies that, by (H$_3$), there exists a unique $x_*\gg 1$ such that, for $\xi_1\in [x_*,\infty)$
\begin{eqnarray}
& &d'(u_+)-b'(u_+) \notag \\
& &> d'(\phi(\xi_1))-\int_{\mathbb{R}^n} f_\beta(y)
b'(\phi(\xi_1-y_1-c\tau)) dy \notag \\
& &\ge d'(\phi(x_*))-\int_{\mathbb{R}^n} f_{\beta}(y)
b'(\phi(x_*-y_1-c\tau)) dy \notag \\
&&>0. \label{01-18}
\end{eqnarray}

\begin{theorem}[Stability of planar traveling waves with time-delay]  Under assumptions  $(H_1)$-$(H_3)$ and $(J_1)$-$(J_2)$,
for a given traveling wave $\phi(x_1+ct)$ of the equation
\eqref{1.1} with $c\ge c_*$ and $\phi(\pm\infty)=u_\pm$, if the
initial data $u_0(s,x)$ is bounded in $[u_-,u_+]$ and $u_0-\phi\in
C([-\tau,0]; H^{m}_w(\mathbb{R}^n)\cap L^1_w(\mathbb{R}^n))$ and
$\partial_s (u_0-\phi) \in L^1([-\tau,0]; H^{m}_w(\mathbb{R}^n)\cap
L^1_w(\mathbb{R}^n))$ with $m>\frac{n}{2}$, then the solution of
\eqref{1.1}  uniquely exists and satisfies:
\begin{enumerate}
\item[$\bullet$] When $c>c_*$,  the solution $u(t,x)$ converges to the noncritical
planar traveling wave $\phi(x_1+ct)$ exponentially
\begin{equation}
\sup_{x\in \mathbb{R}^n}|u(t,x)-\phi(x_1+ct)|\le
C(1+t)^{-\frac{n}{\alpha}}e^{-\mu_\tau t}, \ \ t> 0, \label{2.12}
\end{equation}
where
\begin{equation}\label{mu-tau}
0<\mu_\tau <  \min\{d'(u_+)-b'(u_+), \  \varepsilon_1[G_c(\lambda_*)-H_c(\lambda_*)]\},
\end{equation}
 and $\varepsilon_1=\varepsilon_1(\tau)$ such that $0<\varepsilon_1<1$  for $\tau>0$, and $\varepsilon_1=\varepsilon_1(\tau)\to 0^+$ as $\tau\to+\infty$;

\item[$\bullet$] When $c=c_*$,  the solution $u(t,x)$ converges to the critical planar
traveling wave $\phi(x_1+c_*t)$ algebraically
\begin{equation}
\sup_{x\in \mathbb{R}^n}|u(t,x)-\phi(x_1+c_*t)|\le
C(1+t)^{-\frac{n}{\alpha}}, \ \ t> 0. \label{2.12-2}
\end{equation}
\end{enumerate}
 \label{thm1}
\end{theorem}

However, when the time-delay $\tau=0$,  then we have the following stronger stability for
the traveling waves but with a weaker condition on initial perturbation.

\begin{theorem}[Stability of planar traveling waves without time-delay]  Under assumptions  $(H_1)$-$(H_3)$ and $(J_1)$-$(J_2)$,
for a given traveling wave $\phi(x_1+ct)$
of the equation \eqref{1.1-2} with $c\ge c_*$ and
$\phi(\pm\infty)=u_\pm$, if the initial data $u_0(x)$ is bounded in $[u_-,u_+]$
and $u_0-\phi\in  H^{m}_w(\mathbb{R}^n)\cap L^1_w(\mathbb{R}^n)$
with $m>\frac{n}{2}$, then the solution of \eqref{1.1-2}  uniquely exists and
satisfies:
\begin{enumerate}
\item[$\bullet$] When $c>c_*$,  the solution $u(t,x)$ converges to the noncritical
planar traveling wave $\phi(x_1+ct)$ exponentially
\begin{equation}
\sup_{x\in \mathbb{R}^n}|u(t,x)-\phi(x_1+ct)|\le
C(1+t)^{-\frac{n}{\alpha}}e^{-\mu_0 t}, \ \ t> 0, \label{2.12-new-1}
\end{equation}
where
\begin{equation}\label{mu-0}
0<\mu_0 <  \min\{d'(u_+)-b'(u_+), \  G_c(\lambda_*)-H_c(\lambda_*)\};
\end{equation}

\item[$\bullet$] When $c=c_*$,  the solution $u(t,x)$ converges to the critical planar
traveling wave $\phi(x_1+c_*t)$ algebraically
\begin{equation}
\sup_{x\in \mathbb{R}^n}|u(t,x)-\phi(x_1+c_*t)|\le
C(1+t)^{-\frac{n}{\alpha}}, \ \ t> 0. \label{2.12-new-2}
\end{equation}
\end{enumerate}
 \label{thm2}
\end{theorem}

\begin{remark}

\

\begin{enumerate}
\item
Comparing Theorem \ref{thm1} with time-delay and Theorem \ref{thm2} without time-delay, we realize that, the sufficient condition on the initial perturbation around the wave in the case with time-delay is stronger than the case without time-delay, but the convergence rate to the noncritical waves $\phi(x_1+ct)$ for $c>c_*$ in the case with time-delay is weaker than the case without time-delay, see \eqref{mu-tau} for $\mu_\tau \le \varepsilon_1 [G_c(\lambda_*)-H_c(\lambda_*)]< G_c(\lambda_*)-H_c(\lambda_*)$, and \eqref{mu-0} for $\mu_0 \le G_c(\lambda_*)-H_c(\lambda_*)$, and $\varepsilon_1=\varepsilon_1(\tau)\to 0^+$ as $\tau\to +\infty$.  This means, the time-delay  $\tau>0$ effects the stability of traveling waves a lot, not only the higher requirement for the initial perturbation, but also
 the slower convergence rate  for the solution to the noncritical traveling waves.

 \item  The convergence rates showed both in Theorem \ref{thm1} and Theorem \ref{thm2} are explicit and optimal, particularly, the algebraic decay rates for the solution converging to the critical waves. Actually, all of them are derived from the linearized equations.
\end{enumerate}
\end{remark}

\section{Linearized Nonlocal dispersion Equations}\label{reults}

In this section, we will derive the solution formulas for  the
linearized nonlocal dispersion equations with or without
time-delay, as well as their optimal decay rates, which will play a
key role in the stability proof in section 5.

Now let us introduce the solution formula for linear delayed ODEs  \cite{KIK} and the asymptotic behaviors of the solutions \cite{Mei-Wang}.

\begin{lemma}[\cite{KIK}]\label{lemKIK}
Let $z(t)$ be the solution to the following linear time-delayed ODE
with time-delay $\tau>0$
\begin{equation}\label{p1}
\begin{cases}
\displaystyle\frac{d}{dt}z(t)+k_1 z(t) =k_2 z(t-\tau) \\
z(s)=z_0(s), \ \ \ s\in[-\tau,0].
\end{cases}
\end{equation}
Then
\begin{equation}
z(t)=e^{-k_1(t+\tau)} e^{{\bar k_2}t}_\tau z_0(-\tau)
+\int^0_{-\tau} e^{-k_1(t-s)}e^{{\bar k_2}(t-\tau-s)}_\tau
[z_0'(s)+k_1 z_0(s)]ds , \label{p2}
\end{equation}
where
\begin{equation}
{\bar k_2}:=k_2 e^{k_1 \tau}, \label{k_2}
\end{equation}
and $e^{{\bar k_2}t}_\tau$ is the so-called {\tt delayed exponential
function} in the form
\begin{equation}
e^{{\bar k_2} t}_\tau =\begin{cases}
0, & -\infty<t<-\tau, \\
1, & -\tau\le t<0, \\
1+\frac{{\bar k_2} t}{1!}, & 0\le t<\tau, \\
1+\frac{{\bar k_2} t}{1!}+\frac{{\bar k_2}^2(t-\tau)^2}{2!}, & \tau\le t<2\tau, \\
\vdots & \vdots \\
1+\frac{{\bar k_2} t}{1!}+\frac{{\bar k_2}^2(t-\tau)^2}{2!}+\cdots + \frac{{\bar k_2}^m[t-(m-1)\tau]^m}{m!}, & (m-1)\tau\le t<m\tau, \\
\vdots & \vdots
\end{cases}
\label{p3}
\end{equation}
and $e^{\bar{k}_2t}_\tau$ is the fundamental solution to
\begin{equation}\label{p4}
\begin{cases}
\displaystyle\frac{d}{dt}z(t) ={\bar k_2} z(t-\tau) \\
z(s)\equiv 1, \ \ \ s\in[-\tau,0].
\end{cases}
\end{equation}
\end{lemma}

\begin{lemma}[\cite{Mei-Wang}]\label{lemma3} Let $k_1\ge 0$ and $k_2\ge 0$. Then the solution $z(t)$  to \eqref{p1} (or equivalently \eqref{p2}) satisfies
\begin{equation}
|z(t)|\le C_0 e^{-k_1 t} e^{{\bar k_2}t}_\tau, \label{12-26-1}
\end{equation}
where
\begin{equation}
C_0:= e^{-k_1\tau}|z_0(-\tau)| + \int^0_{-\tau} e^{k_1 s}|z'_0(s)+k_1 z_0(s)| ds, \label{12-26-2}
\end{equation}
and the fundamental solution $e^{{\bar k_2} t}_\tau $ with ${\bar k_2}>0$ to \eqref{p4} satisfies
\begin{equation}
e^{{\bar k_2} t}_\tau \le C(1+t)^{-\gamma} e^{{\bar k_2} t},   \ \ t>0,
\label{p5}
\end{equation}
for arbitrary number $\gamma>0$.

Furthermore, when $k_1\ge k_2\ge 0$,  there exists a  constant $\varepsilon_1=\varepsilon_1(\tau)$ with $0< \varepsilon_1 <1$
for $\tau>0$, and
$\varepsilon_1=1$ for $\tau=0$, and $\varepsilon_1=\varepsilon_1(\tau)\to 0^+$ as $\tau\to +\infty$,
such that
\begin{equation}
e^{-k_1 t}e^{{\bar k_2}t}_\tau \le C e^{-\varepsilon_1(k_1-k_2) t}, \ \ t>0,
\label{12-26-3}
\end{equation}
and the solution $z(t)$ to \eqref{p1} satisfies
\begin{equation}
|z(t)|\le C e^{-\varepsilon_1(k_1-k_2) t}, \ \ t>0.
\label{12-26-4}
\end{equation}
\end{lemma}

Now, we consider the following linearized nonlocal time-delayed
dispersion equation (which will be derived in section 5 for
the proof of stability of traveling wavefronts)
\begin{eqnarray} \label{p12}
\begin{cases}
\displaystyle \frac{\partial {v}}{\partial t}
-\int_{\mathbb{R}^n}J(y)e^{-\lambda_\ast y_1}v(t,x-y)dy +c_1
v  \vspace{2mm} \\
\qquad\qquad \displaystyle = c_2 \int_{\mathbb{R}^n} f_\beta(y)
e^{-\lambda_*(y_1+c\tau)}v(t-\tau,x-y) dy ,
\label{2010-13} \vspace{2mm} \\
v(s,x)=v_0(s,x), \ \ s\in [-\tau,0],  \ x\in \mathbb{R}^n
\end{cases}
\end{eqnarray}
for some given constant coefficients $c$, $c_1$ and $c_2$, where
$c\ge c_*$ is the wave speed.

We are going to derive its solution formula as well as the
asymptotic behavior of the solution. By taking Fourier transform to
\eqref{p12},  and noting that,
\begin{eqnarray}
& &\mathcal{F}\Big[\int_{\mathbb{R}^n} J(y)e^{-\lambda_* y_1} v(t,x-y) dy\Big](t,\eta) \notag \\
& &=\int_{\mathbb{R}^n} e^{-\mbox{i}x\cdot \eta}\Big(\int_{\mathbb{R}^n} J(y)e^{-\lambda_* y_1} v(t,x-y) dy\Big) dx \nonumber \\
& &=\int_{\mathbb{R}^n} J(y)e^{-\lambda_* y_1} \Big(\int_{\mathbb{R}^n} e^{-\mbox{i}x\cdot \eta} v(t,x-y)dx \Big) dy \notag \\
& &=\int_{\mathbb{R}^n} J(y)e^{-\lambda_* y_1} \Big(\int_{\mathbb{R}^n} e^{-\mbox{i}(x+y)\cdot \eta} v(t,x)dx \Big) dy \notag \\
& &=\Big(\int_{\mathbb{R}^n} e^{-\mbox{i}y\cdot
\eta}J(y)e^{-\lambda_* y_1}dy\Big) \hat{v}(t,\eta), \label{p13-0}
\end{eqnarray}
and
\begin{eqnarray}\label{p15}
& & \mathcal{F}\Big[c_2 \int_{\mathbb{R}^n} f_\beta (y) e^{-\lambda_*(y_1+c \tau)} v(t-\tau,x-y) dy\Big](t-\tau,\eta) \notag \\
& &=c_2\int_{\mathbb{R}^n} e^{-\mbox{i}x\cdot \eta}
\Big(\int_{\mathbb{R}^n} f_\beta (y)
e^{-\lambda_*(y_1+c \tau)} v(t-\tau,x-y) dy\Big) dx \notag \\
& &=c_2 \int_{\mathbb{R}^n} f_\beta (y) e^{-\lambda_*(y_1+c \tau)}
\Big( \int_{\mathbb{R}^n} e^{-\mbox{i}x\cdot \eta} v(t-\tau,x-y)
dx\Big) dy
 \notag \\
& &=c_2 \int_{\mathbb{R}^n} f_\beta (y) e^{-\lambda_*(y_1+c \tau)}
\Big( \int_{\mathbb{R}^n} e^{-\mbox{i}(x+y)\cdot \eta}  v(t-\tau,x)
dx\Big) dy
 \notag \\
 & &=c_2 \int_{\mathbb{R}^n} f_\beta (y) e^{-\lambda_*(y_1+c \tau)} e^{-\mbox{i}y\cdot \eta}\Big( \int_{\mathbb{R}^n} e^{-\mbox{i}x\cdot \eta}  v(t-\tau,x) dx\Big) dy \notag \\
 & &=\Big(c_2 \int_{\mathbb{R}^n} f_\beta (y) e^{-\lambda_*(y_1+c \tau)} e^{-\mbox{i}y\cdot \eta}dy\Big) \hat{v}(t-\tau,\eta),
 \end{eqnarray}
we have
\begin{equation} \label{p16}
\frac{d \hat{v}}{dt} + A(\eta)\hat{v}=B(\eta) \hat{v}(t-\tau,\eta),
\ \ \mbox{ with } \hat{v}(s,\eta)=\hat{v}_0(s,\eta), \
s\in[-\tau,0],
\end{equation}
where
\begin{eqnarray}\label{p15-2}
A(\eta):=c_1-\int_{\mathbb{R}^n} J(y)e^{-\lambda_\ast
y_1}e^{-\mbox{i}y\cdot \eta}dy
\end{eqnarray}
and
\begin{equation}
B(\eta):=c_2 \int_{\mathbb{R}^n} f_\beta (y) e^{-\lambda_*(y_1+c
\tau)} e^{-\mbox{i}y\cdot \eta}dy. \label{p15-1}
\end{equation}
By using the formula of the delayed ODE  \eqref{p2} in Lemma \ref{lemKIK}, we then solve \eqref{p16} as follows
\begin{eqnarray}
\hat{v}(t,\eta)&=& e^{-A(\eta)(t+\tau)}e^{\mathcal{B}(\eta)t}_\tau \hat{v}_0(-\tau,\eta) \notag \\
&
&+\int^0_{-\tau}e^{-A(\eta)(t-s)}e^{\mathcal{B}(\eta)(t-\tau-s)}_\tau
\Big[\partial_s \hat{v}_0(s,\eta) + A(\eta) \hat{v}_0(s,\eta)\Big]ds,
\label{p17}
\end{eqnarray}
where
\begin{equation}
\mathcal{B(\eta)}:=B(\eta) e^{A(\eta)\tau}. \label{12-26-7}
\end{equation}
Then, by taking the inverse Fourier transform to \eqref{p17}, we get
\begin{eqnarray}
v(t,x)&=& \frac{1}{(2\pi)^{n}}\int_{\mathbb{R}^n} e^{\mbox{i}x\cdot\eta} e^{-A(\eta)(t+\tau)}
e^{\mathcal{B(\eta)}t}_\tau \hat{v}_0(-\tau,\eta) d\eta
\notag \\
& &+\int^0_{-\tau} \frac{1}{(2\pi)^{n}}\int_{\mathbb{R}^n} e^{\mbox{i}x\cdot\eta} e^{-A(\eta)(t-s)}e^{\mathcal{B(\eta)}(t-\tau-s)}_\tau \notag \\
& & \ \ \times \Big[\partial_s\hat{v}_0(s,\eta) + A(\eta) \hat{v}_0(s,\eta)\Big] d\eta
ds,  \label{p18-new}
\end{eqnarray}
and its derivatives
\begin{eqnarray}
\partial^k_{x_j}v(t,x)&=& \frac{1}{(2\pi)^{n}}\int_{\mathbb{R}^n} e^{\mbox{i}x\cdot\eta}(\mbox{i}\eta_j)^k
e^{-A(\eta)(t+\tau)}e^{\mathcal{B(\eta)}t}_\tau \hat{v}_0(-\tau,\eta) d\eta \notag \\
& &+\int^0_{-\tau} \frac{1}{(2\pi)^{n}}\int_{\mathbb{R}^n} e^{\mbox{i}x\cdot\eta}(\mbox{i}\eta_j)^k e^{-A(\eta)(t-s)}e^{\mathcal{B(\eta)}(t-\tau-s)}_\tau  \notag \\
& & \ \ \ \ \ \ \ \ \ \ \ \times \Big[\partial_s\hat{v}_0(s,\eta)+
A(\eta)\hat{v}_0(s,\eta)\Big]d\eta ds \label{p18}
\end{eqnarray}
for $k=0,1,\cdots$ and $ j=1,\cdots, n$.

Now we are going to derive the asymptotic behavior of $v(t,x)$.

\begin{proposition}[Optimal decay rates for $\tau>0$]\label{lemma4}
Suppose that $v_0\in C([-\tau,0]; H^{m+1}(\mathbb{R}^n)\cap
L^1(\mathbb{R}^n))$ and $\partial_s v_0\in L^1([-\tau,0]; H^{m}(\mathbb{R}^n)\cap
L^1(\mathbb{R}^n))$ for $m\ge 0$, and let
\begin{equation}\label{p21}
\begin{cases}
\displaystyle \tilde{c}_1:=c_1-\int_{\mathbb{R}^n}J(y)e^{-\lambda_\ast y_1}dy,  \vspace{2mm}\\
\displaystyle c_3:=c_2 \int_{\mathbb{R}^n} f_\beta (y)
e^{-\lambda_*(y_1+c\tau)}dy>0.
\end{cases}
\end{equation}
If $\tilde{c}_1\geq c_3$, then there exists a constant $\varepsilon_1=\varepsilon_1(\tau)$ as showed in \eqref{12-26-3} satisfying
$0<\varepsilon_1<1$ for $\tau>0$,  such that  the solution of the linearized
equation \eqref{p12} satisfies
\begin{equation}
\| \partial^k_{x_j}v(t)\|_{L^2(\mathbb{R}^n)}\le C \mathcal{E}^k_{v_0}
t^{-\frac{n+2k}{2\alpha}}e^{-\varepsilon_1(\tilde{c}_1-c_3)t}, \
t>0, \label{4.14}
\end{equation}
for $k=0,1, \cdots, [m]$ and $j=1,\cdots, n$,
where
\begin{eqnarray}
\mathcal{E}^k_{v_0}:&=& \|v_0(-\tau)\|_{L^1(\mathbb{R}^n)}+\|v_0(-\tau)\|_{H^{k}(\mathbb{R}^n)} \notag \\
&
&+\int_{-\tau}^{0}[\|(v'_{0s},v_0)(s)\|_{L^{1}(\mathbb{R}^n)}+\|(v'_{0s},v_0)(s)\|_{H^{k}(\mathbb{R}^n)}]ds.
\label{p20-1-1}
\end{eqnarray}
Furthermore, if $m>\frac{n}{2}$, then
\begin{equation}
\| v(t)\|_{L^\infty(\mathbb{R}^n)} \leq
C \mathcal{E}^m_{v_0}
t^{-\frac{n}{\alpha}}e^{-\varepsilon_1(\tilde{c}_1-c_3)t}, \ \ t>0.
\label{p20}
\end{equation}
Particularly, when $\tilde{c}_1=c_3$, then
\begin{equation}
\| v(t)\|_{L^\infty(\mathbb{R}^n)} \leq
C \mathcal{E}^m_{v_0}
t^{-\frac{n}{\alpha}}, \ \ t>0.
\label{p20-nnew}
\end{equation}
\end{proposition}
{\bf Proof}. Let
\begin{eqnarray}
I_1(t,\eta):&=& (\mbox{i}\eta_j)^{k}e^{-A(\eta)(t+\tau)}e^{\mathcal{B(\eta)}t}_\tau \hat{v}_0(-\tau,\eta), \label{p22} \\
I_2(t-s,\eta):&=& (\mbox{i}\eta_j)^{k} e^{-A(\eta)(t-s)}e^{\mathcal{B(\eta)}(t-\tau-s)}_\tau
\Big[\partial_s\hat{v}_0(s,\eta)+ A(\eta)\hat{v}_0(s,\eta)\Big].
\label{p23}
\end{eqnarray}
Then, \eqref{p18} is reduced to
\begin{equation}
\partial^{k}_{x_j} v(t,x)=\mathcal{F}^{-1}[I_1](t,x) + \int^0_{-\tau} \mathcal{F}^{-1}[I_2](t-s,x) ds.
\label{p24}
\end{equation}
So, by using Parseval's equality, we have
\begin{eqnarray}
\|\partial^{k}_{x_j} v(t)\|_{L^2(\mathbb{R}^n)}&\le
&\|\mathcal{F}^{-1}[I_1](t)\|_{L^2(\mathbb{R}^n)} + \int^0_{-\tau}
\|\mathcal{F}^{-1}[I_2](t-s)\|_{L^2(\mathbb{R}^n)} ds
\notag \\
&=&\|I_1(t)\|_{L^2(\mathbb{R}^n)}+\int^0_{-\tau}
\|I_2(t-s)\|_{L^2(\mathbb{R}^n)} ds. \label{p25}
\end{eqnarray}
\begin{eqnarray}\label{p26-1}
|e^{-A(\eta)t}|&=&e^{-c_1t} \Big|\exp\Big(t\int_{\mathbb{R}^n}
J(y)e^{-\lambda_\ast y_1}e^{-\mbox{i}y\cdot\eta}dy\Big)\Big|\nonumber\\
&=&e^{-{c}_1t}\exp\Big(t\int_{\mathbb{R}^n}
J(y)e^{-\lambda_\ast y_1} \cos (y\cdot \eta) dy\Big)\nonumber\\
& =& e^{-\tilde{c}_1t}\exp\Big(-t\int_{\mathbb{R}^n}
J(y)e^{-\lambda_\ast y_1} (1- \cos (y\cdot \eta)) dy\Big)\nonumber\\
&=: &e^{-k_1t},  \ \ \ \ \ \mbox{ with }
k_1:=\tilde{c}_1+\int_{\mathbb{R}^n} J(y)e^{-\lambda_\ast y_1} (1-
\cos (y\cdot \eta)) dy ,
\end{eqnarray}
Note that, using \eqref{p15-2}, \eqref{p15-1},  and the facts
$\frac{e^x+e^{-x}}{2}\ge 1$ for all $x\in \mathbb{R}$, and
$\int_{\mathbb{R}^n}J(y)\sin(y\cdot \eta) dy =0$ because $J(y)$ is
even and $\sin(y\cdot \eta)$ is odd, and
$\int_{\mathbb{R}^n}J(y)dy=1$, we have
\begin{eqnarray} \label{p26}
&&\exp\Big(-t\int_{\mathbb{R}^n} J(y)e^{-\lambda_\ast y_1} (1- \cos
(y\cdot \eta)) dy\Big)\nonumber\\
&& = \exp\Big(-t\int_{\mathbb{R}^n} J(y)\frac{e^{-\lambda_\ast y_1}
+ e^{\lambda_\ast y_1} }{2} (1- \cos
(y\cdot \eta)) dy\Big)\nonumber\\
&& \leq \exp\Big(-t\int_{\mathbb{R}^n}
J(y)(1- \cos (y\cdot \eta)) dy\Big)\nonumber\\
&&= \exp\Big(-t\int_{\mathbb{R}^n}
J(y)[1- [\cos (y\cdot \eta)+\mbox{i}\sin (y\cdot\eta)]] dy\Big)\nonumber\\
&&= e^{(\hat{J}(\eta)-1)t}
\end{eqnarray}
and
\begin{equation}
|B(\eta)|\le c_2 \int_{\mathbb{R}^n} f_\beta (y)
e^{-\lambda_*(y_1+c\tau)}dy=c_3=:k_2, \label{p27}
\end{equation}
 and
\begin{equation}
|\mathcal{B}(\eta)|=|B(\eta)e^{A(\eta) \tau}|\le c_3 e^{k_1\tau}=k_2
e^{k_1 \tau}=:{\bar k}_2, \label{p27-2}
\end{equation}
and further
\begin{equation}
|e^{\mathcal{B(\eta)}t}_\tau|\le e^{\bar k_2t}_\tau. \label{p28}
\end{equation}
If $\tilde{c}_1\ge c_3$, from (J$_2$), namely,
$1-\hat{J}(\eta)=\mathcal{K}|\eta|^\alpha-o(|\eta|^\alpha)>0\
\mbox{as}\ \eta\rightarrow0 $, then $k_1=\tilde{c}_1
+1-\hat{J}(\eta)\geq c_3=k_2$. Using \eqref{p26-1}, \eqref{p26},
\eqref{p28} and \eqref{12-26-3} in Lemma \ref{lemma3}, we obtain
\begin{eqnarray}
\|I_1(t)\|^2_{L^2(\mathbb{R}^n)}&=&\int_{\mathbb{R}^n}
|e^{-A(\eta)(t+\tau)}e^{\mathcal{B}(\eta)t}_\tau
\hat{v}_0(-\tau,\eta)|^2|\eta_j|^{2k} d\eta\nonumber\\
&&\leq
C\int_{\mathbb{R}^n}(e^{-k_1(t+\tau)}e^{\bar{k}_2t}_\tau)^2|\hat{v}_0(-\tau,\eta)|^2|\eta_j|^{2k}
d\eta\nonumber\\
& &\le C\int_{\mathbb{R}^n}(e^{-\varepsilon_1(k_1-k_2)t})^2|\hat{v}_0(-\tau,\eta)|^2|\eta_j|^{2k}
d\eta\nonumber\\
 &&= Ce^{-2\varepsilon_1(\tilde{c}_1-c_3)t}
\int_{\mathbb{R}^n}e^{-2\varepsilon_1(1-\hat{J}(\eta))t}|\hat{v}_0(-\tau,\eta)|^2|\eta_j|^{2k}
d\eta. \label{p29}
\end{eqnarray}
Again from (J$_2$),   there exist
some numbers $0<\mathcal{K}_1<\mathcal{K}$, $0<\delta<1$ and  $\tilde{a}>0$, such that
\begin{eqnarray}\label{4.11}
\begin{cases}
\mathcal{K}_1|\eta|^{\alpha}\le 1-\hat{J}(\eta) \leq \mathcal{K}|\eta|^{\alpha}, & \mbox{as}\
|\eta|\leq \tilde{a},\\
\delta:=\mathcal{K}_1{\tilde a}^{\alpha}\le 1-\hat{J}(\eta) \leq \mathcal{K}|\eta|^{\alpha}, & \mbox{as}\ |\eta|\geq \tilde{a}.
\end{cases}
\end{eqnarray}
Therefore, we have
\begin{eqnarray}\label{4.12}
&&\int_{\mathbb{R}^n}e^{-2\varepsilon_1 (1-\hat{J}(\eta))t}|\hat{v}_0(-\tau,\eta)|^2|\eta_j|^{2k}
d\eta\nonumber\\
&&=\int_{|\eta|\leq
\tilde{a}}e^{-2\varepsilon_1(1-\hat{J}(\eta))t}|\hat{v}_0(-\tau,\eta)|^2|\eta_j|^{2k}
d\eta+\int_{|\eta|\geq
\tilde{a}}e^{-2\varepsilon_1(1-\hat{J}(\eta))t}|\hat{v}_0(-\tau,\eta)|^2|\eta_j|^{2k}
d\eta\nonumber\\
&&\leq\int_{|\eta|\leq
\tilde{a}}e^{-2\varepsilon_1\mathcal{K}_1|\eta|^{\alpha}t}|\hat{v}_0(-\tau,\eta)|^2|\eta_j|^{2k}
d\eta+\int_{|\eta|\geq \tilde{a}}e^{-2\varepsilon_1\delta
t}|\hat{v}_0(-\tau,\eta)|^2|\eta_j|^{2k}
d\eta\nonumber\\
&&\leq \|\hat{v}_0(-\tau)\|_{L^\infty(\mathbb{R}^n)}^2
t^{-\frac{n+2k}{\alpha}}\int_{|\eta|\leq \tilde{a}}e^{-2\varepsilon_1 \mathcal{K}_1|\eta
t^{\frac1\alpha}|^{\alpha}}|\eta_jt^{\frac1\alpha}|^{2k}
d(\eta t^{\frac{1}{\alpha}}) \notag \\
& & \ \ \ +e^{-2\varepsilon_1 \delta t}\int_{|\eta|\geq
\tilde{a}}|\hat{v}_0(-\tau,\eta)|^2|\eta_j|^{2k}
d\eta\nonumber\\
&&\leq C(\|v_0(-\tau)\|_{L^1(\mathbb{R}^n)}^2+\|v_0(-\tau)\|_{H^k(\mathbb{R}^n)}^2)
t^{-\frac{n+2k}{\alpha}}.
\end{eqnarray}
Substitute \eqref{4.12} into  \eqref{p29}, we obtain
\begin{eqnarray}\label{4.13}
\|I_1(t)\|_{L^2(\mathbb{R}^n)}\leq
C(\|v_0(-\tau)\|_{L^1(\mathbb{R}^n)}+\|v_0(-\tau)\|_{H^k(\mathbb{R}^n)})t^{-\frac{n+2k}{2\alpha}}e^{-\varepsilon_1(\tilde{c}_1-c_3)t}
\end{eqnarray}

 Thus, in a similar way, we can also prove
\begin{eqnarray}\label{p32}
& &\|I_2(t-s)\|_{L^2(\mathbb{R}^n)}\nonumber\\
&&=\bigg(\int_{\mathbb{R}^n}
|e^{-A(\eta)(t-s)}e^{\mathcal{B(\eta)}(t-\tau-s)}_\tau|^2
\Big|\partial_s \hat{v}_0(s,\eta) + A(\eta)\hat{v}_0(s,\eta)\Big|^2\cdot|\eta_j|^{2k}d\eta\bigg)^{\frac12}\notag \\
& &\leq Ce^{-\varepsilon_1(\tilde{c}_1-c_3)t}
\bigg(\int_{\mathbb{R}^n}e^{-2\varepsilon_1(1-\hat{J}(\eta))t}
\Big(|\eta|^{2k}|\partial_s\hat{v}_0(s,\eta)|+|\eta|^{2k}|\hat{v}_0(s,\eta)|^2 \Big)d\eta\bigg)^{\frac12} \notag \\
& &\leq
Ct^{-\frac{n+2k}{2\alpha}}e^{-\varepsilon_1(\tilde{c}_1-c_3)t}
\Big(\|(\partial_s v_0,v_0)(s)\|_{L^{1}(\mathbb{R}^n)}+\|(\partial_s
v_0,v_0)(s)\|_{H^{k}(\mathbb{R}^n)}\Big).
\end{eqnarray}
Substituting \eqref{4.13} and \eqref{p32} to \eqref{p25}, we
immediately obtain \eqref{4.14}.

Similarly, we can prove \eqref{p20}. We omit the details. Thus, we
complete the proof of Proposition \ref{lemma4}.  $\square$

\medskip

For $\tau=0$, the equation \eqref{p12} is reduced to
\begin{eqnarray} \label{p12-new2}
\begin{cases}
\displaystyle \frac{\partial {v}}{\partial t}+c\frac{\partial v}{\partial x_1}
-\int_{\mathbb{R}^n}J(y)e^{-\lambda_\ast y_1}v(t,x-y)dy +c_1
v  \vspace{2mm} \\
\qquad\qquad \displaystyle = c_2 \int_{\mathbb{R}^n} f_\beta(y)
e^{-\lambda_*(y_1+c\tau)}v(t,x-y-c\tau{\bf e}_1) dy ,
\label{2010-13-new} \vspace{2mm} \\
v(s,x)=v_0(x), \ \  \ x\in \mathbb{R}^n.
\end{cases}
\end{eqnarray}
Taking Fourier transform to \eqref{p12-new2}, as showed in \eqref{p16}, we have
\begin{equation} \label{p16-new2}
\frac{d \hat{v}}{dt} =[B(\eta)- A(\eta)]\hat{v},
\ \ \mbox{ with } \hat{v}(0,\eta)=\hat{v}_0(\eta),
\end{equation}
where $A(\eta)$ and $B(\eta)$ are given in \eqref{p15-2} and \eqref{p15-1} with $\tau=0$, respectively. Integrating \eqref{p16-new2} yields
\[
\hat{v}(t,\eta)=e^{-[A(\eta)-B(\eta)]t}\hat{v}_0(\eta).
\]
Taking the inverse Fourier transform, we get the solution formula
\[
v(t,x)=\frac{1}{(2\pi)^n}\int_{\mathbb{R}^n} e^{\mbox{i}x\cdot \eta} e^{-[A(\eta)-B(\eta)]t}\hat{v}_0(\eta) d\eta.
\]
Then, a similar analysis as showed before can derive the optimal decay of the solution in the case without time-delay as follows. The detail of proof is omitted.

\begin{proposition}[Optimal decay rates for $\tau=0$]\label{lemma5}
Suppose that $v_0\in  H^m(\mathbb{R}^n)\cap
L^1(\mathbb{R}^n))$ for $m\ge 0$, then  the solution of the linearized
equation \eqref{p12-new2} satisfies
\begin{equation}
\| \partial^k_{x_j}v(t)\|_{L^2(\mathbb{R}^n)}\le C (\|v_0\|_{L^1(\mathbb{R}^n)}+\|v_0\|_{H^k(\mathbb{R}^n)})
t^{-\frac{n+2k}{2\alpha}}e^{-(\tilde{c}_1-c_3)t}, \
t>0, \label{4.14-new2}
\end{equation}
for $k=0,1, \cdots, [m]$ and $j=1,\cdots, n$, where the positive constants $\tilde{c}_1$ and $c_3$ are defined in \eqref{p21} for $\tau=0$.

Furthermore, if $m>\frac{n}{2}$, then
\begin{equation}
\| v(t)\|_{L^\infty(\mathbb{R}^n)} \leq
C (\|v_0\|_{L^1(\mathbb{R}^n)}+\|v_0\|_{H^k(\mathbb{R}^n)})
t^{-\frac{n}{\alpha}}e^{-(\tilde{c}_1-c_3)t}, \ \ t>0.
\label{p20-new2}
\end{equation}
Particularly, when $\tilde{c}_1=c_3$, then
\begin{equation}
\| v(t)\|_{L^\infty(\mathbb{R}^n)} \leq
C (\|v_0\|_{L^1(\mathbb{R}^n)}+\|v_0\|_{H^k(\mathbb{R}^n)})
t^{-\frac{n}{\alpha}}, \ \ t>0.
\label{p20-nmmnew}
\end{equation}
\end{proposition}

\section{Global Existence  and Comparison principle}

In this section, we prove the global existence and uniqueness of the solution for the Cauchy problem
to the nonlinear equation with nonlocal dispersion \eqref{1.1}, and then establish the comparison principle
in $n$-D case by a different proof approach to the previous work \cite{Chen,Coville-Dupaigne-2}.

\begin{proposition}[Existence and Uniqueness]\label{lemma6}
Let $u_0(s,x)\in C([-\tau,0]; C( \mathbb{R}^n))$ with
$0=u_-\le u_0(s,x)\le u_+$ for $(s,x)\in [-\tau,0]\times \mathbb{R}^n$, then the  solution to \eqref{1.1}
 uniquely and globally exists, and satisfies that
$u\in C^1([0,\infty); C( \mathbb{R}^n))$, and
$ u_-\le u(t,x) \le u_+$ for $ (t,x)\in \mathbb{R}_+\times \mathbb{R}^n)$.
\end{proposition}
\noindent {\bf Proof}. Multiplying \eqref{1.1} by $e^{\eta_0 t}$ and integrating it over $[0,t]$ with respect to $t$, where
$\eta_0>0$ will be technically selected in \eqref{01-14} below, we  then express
\eqref{1.1}  in the integral form
\begin{eqnarray}\label{01-15}
u(t,x)&=&e^{-\eta_0 t}u(0,x)+\int_{0}^{t}e^{-\eta_0
(t-s)}\bigg[\int_{\mathbb{R}^n}J(x-y)u(s,y)dy+(\eta_0-1)
u(s,x)\nonumber\\
&& \ \ \ -d(u(s,x))+\int_{\mathbb{R}^n}f_\beta(y)b(u(s-\tau,x-y))dy\bigg]ds.
\end{eqnarray}
Let us define the solution space as, for any $T\in[0,\infty]$,
\begin{align}
\mathfrak{B}=\Big\{ u(t,x)|& u(t,x)\in C([0,T]\times \mathbb{R}^n) \mbox{ with }
u_-\le u \le u_+, \notag \\
 & \qquad u(s,x)=u_0(s,x), (s,x)\in [-\tau,0] \times \mathbb{R}^n\Big\}, \label{01-13}
\end{align}
with the norm
\begin{equation}
\|u\|_{\mathfrak{B}}=\sup_{t\in[0,T]}e^{-\eta_0
t}\|u(t)\|_{L^\infty(\mathbb{R}^n)} ,
\end{equation}
where
\begin{equation}
\eta_0: =1+\eta_1 + \eta_2, \ \ \ \
\eta_1:=\max_{u\in[u_-,u_+]}|d'(u)|, \ \ \ \
\eta_2:=\max_{u\in[u_-,u_+]}|b'(u)|.
 \label{01-14}
\end{equation}
Clearly,   $\mathfrak{B}$ is a Banach space.

Define an operator $\mathcal{P}$ on $\mathfrak{B}$ by
\begin{eqnarray}
\mathcal{P}(u)(t,x):&=&
e^{-\eta_0 t}u_0(0,x)+\int_{0}^{t}e^{-\eta_0
(t-s)}\bigg[\int_{\mathbb{R}^n}J(x-y)u(s,y)dy+(\eta_0-1)
u(s,x)\nonumber\\
& & \ \
-d(u(s,x))+\int_{\mathbb{R}^n}f_\beta(y)b(u(s-\tau,x-y))dy\bigg]ds, \   \ \  \mbox{ for }\
\ 0\le t\le T, \label{01-32}
\end{eqnarray}
and
\begin{equation}
\mathcal{P}(u)(s,x):=
u_0(s,x), \ \ \ \mbox{for}\ \ s\in[-\tau,0].
\end{equation}

Now we are going to prove that $\mathcal{P}$ is a contracting operator from $\mathfrak{B}$ to $\mathfrak{B}$.

Firstly, we prove that $\mathcal{P}: \mathfrak{B} \to \mathfrak{B}$. In fact, if $u\in \mathfrak{B}$, from (H$_1$)-(H$_3$), namely, $0=d(0)\le d(u)\le d(u_+)$, $0=b(0)\le b(u)\le b(u_+)$, and $d(u_+)=b(u_+)$,
and using the facts $\int_{\mathbb{R}^n}J(x-y)dy=1$, $\int_{\mathbb{R}^n} f_\beta(y)dy=1$, and
\begin{equation}
g(u):=(\eta_0-1)u-d(u) \mbox{ is increasing},
\label{g}
\end{equation}
which implies $g(u_+)\ge g(u)\ge g(0)=0$ for $u\in [u_-,u_+]$,
then we have
\begin{eqnarray}\label{4.1}
0=u_-\leq \mathcal{P}(u)&\leq& e^{-\eta_0 t}u_+ +\int_{0}^{t}e^{-\eta_0
(t-s)}\bigg[\int_{\mathbb{R}^n}J(x-y)u_+dy  \notag \\
& & \ \ \ \ +(\eta_0-1)
u_+ -d(u_+) +\int_{\mathbb{R}^n}f_\beta(y)b(u_+)dy\bigg]ds \notag \\
&=&e^{-\eta_0 t}u_+ +\int_{0}^{t}e^{-\eta_0
(t-s)}[\eta_0 u_+ - d(u_+)+b(u_+)]ds \notag \\
&=& u_+.
\end{eqnarray}
This plus the continuity of $\mathcal{P}(u)$ based on the continuity of $u$ proves $\mathcal{P}(u)\in\mathfrak{B}$, namely, $\mathcal{P}$ maps from $\mathfrak{B}$ to $\mathfrak{B}$.

Secondly, we prove that $\mathcal{P}$ is contracting. In fact, let $u_1, \ u_2 \in \mathfrak{B}$, and $v=u_1-u_2$, then we have
\begin{eqnarray}\label{4.2}
&&\mathcal{P}(u_1)-\mathcal{P}(u_2)\\
&&=\int_{0}^{t}e^{-\eta_0
(t-s)}\bigg[\int_{\mathbb{R}^n}J(x-y)v(s,y)dy+(\eta_0-1) v(s,x)
-[d(u_1(s,x))-d(u_2(s,x))]\nonumber\\
&&\qquad\quad+\int_{\mathbb{R}^n}f_\beta(y)[b(u_1(s-\tau,x-y))-b(u_2(s-\tau,x-y))]dy\bigg]ds.
\end{eqnarray}
So, we have
\begin{eqnarray}\label{4.3}
|\mathcal{P}(u_1)-\mathcal{P}(u_2)|e^{-\eta_0 t}&\leq & \int_{0}^{t}e^{-2\eta_0
(t-s)}\Big(\eta_0+\max_{u\in[u_-,u_+]}|d'(u)|\Big)\|v\|_{\mathfrak{B}} ds \nonumber\\
&&  +\max_{u\in[u_-,u_+]}|b'(u)|
\begin{cases}
\int_{0}^{t-\tau}e^{-2\eta_0 (t-s)}\|v\|_{\mathfrak{B}}ds, & \mbox{for}\  t\geq \tau \\
0, & \mbox{for}\ 0\leq t\leq \tau
\end{cases} \nonumber \\
&\leq& \frac{1}{2\eta_0}\Big((\eta_0+\eta_1)(1-e^{-2\eta_0
t})+\eta_2(e^{-2\eta_0
\tau}-e^{-2\eta_0 t})\Big)\|v\|_{\mathfrak{B}}\nonumber\\
&\le & \frac{\eta_0+\eta_1+\eta_2}{2\eta_0} \|v\|_{\mathfrak{B}} \notag \\
&=&  \frac{2\eta_0-1}{2\eta_0}\|v\|_{\mathfrak{B}} \notag \\
&=:& \rho \|v\|_{\mathfrak{B}}
\end{eqnarray}
for $0<\rho:=\frac{2\eta_0-1}{2\eta_0}<1$, namely, we prove that the mapping $\mathcal{P}$ is contracting:
\begin{equation}\label{4.4}
\|\mathcal{P}(u_1)-\mathcal{P}(u_2)\|_{\mathfrak{B}}\leq \rho \|u_1-u_2\|_{\mathfrak{B}}<\|u_1-u_2\|_{\mathfrak{B}}.
\end{equation}

Hence, by the Banach fixed-point theorem,  $\mathcal{P}$ has a unique fixed point
$u$ in $\mathfrak{B}$, i.e, the integral equation \eqref{01-15} has a unique classical solution on $[0,T]$ for
any given $T>0$.  Differentiating \eqref{01-15} with respect to $t$, we  get back to the original equation \eqref{1.1}, i.e.,
\begin{equation}\label{01-16}
u_t = J\ast u -u+d(u(t,x))+
\int_{\mathbb{R}^n} f_\beta
(y) b(u(t-\tau,x-y)) dy,
\end{equation}
then we can easily confirm from the right-hand-side of \eqref{01-16} that $u_t \in C([0,T]\times \mathbb{R}^n)$.
 This completes our proof. $\square$

\begin{remark}\label{remark}
From the proof of Proposition \ref{lemma6}, we realize that, when $u_0(s,x)\in C^k([-\tau,0]\times\mathbb{R}^n)$, then the solution of the time-delayed equation \eqref{1.1} holds $u(t,x)\in C^{k+1}([0,\infty);C(\mathbb{R}^n))$; while for the non-delayed equation \eqref{1.1-2} (i.e., $\tau=0$), if $u_0(x)\in C(\mathbb{R}^n)$, then the  solution of  the non-delayed equation \eqref{1.1-2} holds
$u(t,x)\in C^{\infty}([0,\infty);C(\mathbb{R}^n))$. This means that the solution to the nonlocal dispersion equation \eqref{1.1} possesses
a really good regularity in time. However,  the solutions for \eqref{1.1} lack the regularity in space.

\end{remark}

Now we establish two comparison principle for \eqref{1.1}. Although the comparison principle in 1D case were proved in
\cite{Chen, Coville-Dupaigne-2}. Here we give a comparison principle in $n$-D case with much weaker restriction on the initial data. The proof is also new and easy to follow. Different from the previous works \cite{Chen, Coville-Dupaigne-2}, instead of the differential equation \eqref{1.1}, we will work on the integral equation \eqref{01-15}, and sufficiently use the property of contracting operator $\mathcal{P}$.

Let $\bar{u}(t,x)$ be an upper solution to \eqref{1.1}, namely
\begin{equation}
\begin{cases}
\displaystyle \frac{\partial {\bar u}}{\partial t} - J\ast {\bar u} +{\bar u}+d({\bar u}(t,x))\ge
\int_{\mathbb{R}^n} f_\beta
(y) b({\bar u}(t-\tau,x-y)) dy, \\
{\bar u}(s,x) \ge u_0(s,x), \ \ s\in[-\tau,0], \ x\in \mathbb{R}^n,
\end{cases}
\label{1.1-30}
\end{equation}
where its integral form can be written as
\begin{eqnarray}\label{01-15-new-new}
&&{\bar u}(t,x)\ge e^{-\eta_0 t}{\bar u}(0,x)+\int_{0}^{t}e^{-\eta_0
(t-s)}\bigg[\int_{\mathbb{R}^n}J(x-y){\bar u}(s,y)dy+(\eta_0-1)
{\bar u}(x,s)\nonumber\\
&&\qquad\qquad\qquad\qquad-d({\bar u}(s,x))+\int_{\mathbb{R}^n}f_\beta(y)b({\bar u}(s-\tau,x-y))dy\bigg]ds, \ \ \mbox{ for } t>0
\end{eqnarray}
and let $\underline{u}(t,x)$ be an lower solution to \eqref{1.1} satisfying \eqref{1.1-30} or \eqref{01-15-new-new} conversely.
Then we have the following comparison result.

\begin{proposition}[Comparison Principle]\label{lemma7}
Let $\underline{u}(t,x)$ and $\bar{u}(t,x)$ be the classical lower and upper solutions to \eqref{1.1},  with $u_-\le \underline{u}(t,x), \ \bar{u}(t,x)\le u_+$, respectively, and satisfy $0\le \underline{u}(t,x) \le u_+$ and  $0\le \bar{u}(t,x) \le u_+$ for
$(t,x)\in \mathbb{R}_+\times \mathbb{R}^n$.   Then  $\underline{u}(t,x)\leq
\bar{u}(t,x)$ for $(t,x)\in [0,\infty)\times \mathbb{R}^n$.
\end{proposition}
\noindent{\bf Proof}. We need to prove
$\bar{u}(t,x)-\underline{u}(t,x)\ge 0$ for $(t,x)\in
[0,\infty)\times \mathbb{R}^n$, namely, $r(t):=\inf_{x\in
\mathbb{R}^n} v(t,x)\ge 0$, where
$v(t,x):=\bar{u}(t,x)-\underline{u}(t,x)$.

If this is not true, then there exist some constants $\varepsilon>0$
and $T>0$ such that $r(t)>-\varepsilon e^{3\eta_0t}$ for $t\in[0,T)$
and $r(T)=-\varepsilon e^{3\eta_0T}$, where $\eta_0$ given in
\eqref{01-14}.

Since $\underline{u}(t,x)$ and ${\bar u}(t,x)$ are the lower and
upper solutions to \eqref{1.1} and
$\bar{u}(s,x)-\underline{u}(s,x)\geq0$, for $s\in[-\tau,0]$, and using \eqref{01-14} and \eqref{g}, and noting
${\bar u}(t,x)-\underline{u}(t,x)\ge -\varepsilon e^{3\eta_0 T}$ for $(t,x)\in [0,T]\times \mathbb{R}^n$, then
we have, for $0\le t\le T$,
\begin{eqnarray}
&&\bar{u}(t,x)-\underline{u}(t,x)\nonumber\\
&&\geq e^{-\eta_0 t}[\bar{u}(0,x)-\underline{u}(0,x)] \notag \\
& & \ \ \ +\int^t_0 e^{-\eta_0(t-s)}\Big(\int_{\mathbb{R}^n} J(x-y) [\bar{u}(s,y)-\underline{u}(s,y)]dy \notag \\
& & \ \ \  + g(\bar{u}(s,x))-g(\underline{u}(s,x)) \notag \\
& & \ \ \ +\int_{\mathbb{R}^n} f_\beta(y)[b(\bar{u}(s-\tau,x-y))-b(\underline{u}(s-\tau,x-y))]dy\Big)ds \notag \\
&&\geq\int^t_0 e^{-\eta_0(t-s)}\Big(-\varepsilon e^{3\eta_0s}-\max_{\zeta\in[u_-,u_+]}g'(\zeta)\varepsilon e^{3\eta_0s}\Big)ds\nonumber\\
 &&\ \ \ -\max_{u\in[u_-,u_+]}|b'(u)|
\begin{cases}
\int_{\tau}^{t}e^{-\eta_0 (t-s)}\varepsilon e^{3\eta_0(s-\tau)}ds, & \mbox{for}\  t\geq \tau \\
0, & \mbox{for}\ 0\leq t\leq \tau
\end{cases} \nonumber \\
 &&\geq\begin{cases}
-(\eta_0+1)\varepsilon e^{-\eta_0t}\int^t_0 e^{4\eta_0s}ds-\eta_0\varepsilon e^{-3\eta_0\tau}e^{-\eta_0t}\int_{\tau}^{t} e^{4\eta_0s}ds,  &\mbox{for}\  t\geq \tau \\
-(\eta_0+1)\varepsilon e^{-\eta_0t}\int^t_0 e^{4\eta_0s}ds, &
\mbox{for}\ 0\leq t\leq \tau
\end{cases}\nonumber\\
&&\geq-\frac{2\eta_0+1}{4\eta_0}\varepsilon e^{3\eta_0t}.
\end{eqnarray}
Thus, from the assumption we know
\begin{eqnarray}
-\varepsilon e^{3\eta_0T}=\inf_{x\in
\mathbb{R}^n}(\bar{u}(T,x)-\underline{u}(T,x))\geq-\frac{2\eta_0+1}{4\eta_0}\varepsilon
e^{3\eta_0T},
\end{eqnarray}
which is a contradiction for $\eta_0>\frac{1}{2}$. Here, our $\eta_0$ defined in \eqref{01-14} satisfies
$\eta_0>1$. Thus the proof is complete.  $\square$

\section{Global Stability of Planar Traveling Waves}

 The main purpose in this section is to prove
  Theorems \ref{thm1} for all traveling waves including
 the critical traveling waves.

For given traveling wave $\phi(x_1+ct)$ with the speed $c\ge c_*$ and the given initial data $u_-\le u_0(s,x)\le u_+$ for $(s,x)\in [-\tau,0]\times \mathbb{R}^n$,
let us define $U^+_0(s,x)$ and $U^-_0(s,x)$ as
\begin{eqnarray}
U^-_0(s,x):&=& \min\{ \phi(x_1+cs), u_0(s,x)\} \notag \\
U^+_0(s,x):&=&\max\{  \phi(x_1+cs), u_0(s,x)\}
\label{2.15}
\end{eqnarray}
for $(s,x)\in [-\tau,0]\times \mathbb{R}^n$. So,
\[
u_0-\phi=(U^+_0-\phi) + (U^-_0-\phi).
\]

Since $u_0-\phi\in C([-\tau,0];H^{m+1}_w(\mathbb{R}^n)\cap L^1_w(\mathbb{R}^n))$ with $m>\frac{n}{2}$ and $w(x)\ge 1$ (see \eqref{2.8-new}), and noting Sobolev's embedding theorem
$H^m(\mathbb{R}^n)\hookrightarrow C(\mathbb{R}^n)$, we have $u_0-\phi \in C([-\tau,0];C(\mathbb{R}^n))$. On the other hand,
the traveling wave $\phi(x_1+cs)$ is smooth, then we can guarantee $U_0^\pm(s,x) \in C([-\tau,0];C(\mathbb{R}^n))$. Thus,
applying Proposition \ref{lemma6}, we know that the solutions of \eqref{1.1}
with the initial data $U^+_0(s,x)$ and $U^-_0(s,x)$ globally exist,
and denote them by  $U^+(t,x)$ and $U^-(t,x)$, respectively, that is,
\begin{equation}
\begin{cases}
\displaystyle \frac{\partial U^\pm}{\partial t}-J\ast U^{\pm}+U^\pm+d( U^\pm)=
\int_{\mathbb{R}^n} f_\beta (y) b(U^\pm(t-\tau,x-y)) dy,\\
U^\pm(s,x)=U^\pm_0(s,x), \ \ \  x\in \mathbb{R}^n, \ s\in[-\tau,0].
\label{2.16-1}
\end{cases}
\end{equation}
Then the comparison principle (Proposition \ref{lemma7}) further implies
\begin{eqnarray}\label{2.18}
\begin{cases}
u_-\le U^-(t,x)\le u(t,x)\le U^+(t,x)\le u_+
  \\
u_-\le  U^-(t,x)\le \phi(x_1+ct)\le U^+(t,x)\le u_+
\end{cases}
\ \ \mbox{ for }(t,x)\in \mathbb{R}_+\times \mathbb{R}^n.
\end{eqnarray}

In what follows, we are going to complete the proof for the
stability in three steps.

{\bf Step 1. The convergence of $U^+(t,x)$ to $\phi(x_1+ct)$}

Let
\begin{equation}
V(t,x):=U^+(t,x)-\phi(x_1+ct), \ \ \ \
V_0(s,x):=U^+_0(s,x)-\phi(x_1+cs) . \label{2.19}
\end{equation}
It follows from \eqref{2.18} that
\begin{equation}
V(t,x)\ge 0, \qquad \ V_0(s,x)\ge 0. \label{2.20}
\end{equation}
We see from \eqref{1.1} that $V(t,x)$ satisfies (by linearizing it
around 0)
\begin{eqnarray}\label{2.21}
& & \displaystyle \frac{\partial V}{\partial t}
-\int_{\mathbb{R}^n}J(y) V(t,x-y)dy+V +d'(0) V \notag \\
& & \quad -  b'(0) \int_{\mathbb{R}^n} f_\beta(y) V(t-\tau,x-y) dy \notag \\
& &  = -Q_1(t,x)+  \int_{\mathbb{R}^n} f_\beta(y) Q_2(t-\tau, x-y)
dy +[d'(0)-d'(\phi(x_1+ct))]V
\notag \\
& &\quad +   \int_{\mathbb{R}^n} f_\beta(y)
[b'(\phi(x_1-y_1+c(t-\tau))-b'(0)]V(t-\tau,x-y) dy \notag \\
 & &=:  I_1(t,x)+I_2(t-\tau,x)+I_3(t,x)+I_4(t-\tau,x),
\end{eqnarray}
with the initial data
\begin{equation}
 V(s,x)=V_0(s,x), \ s\in [-\tau,0],
 \label{2.21-1}
 \end{equation}
 where
\begin{equation}
 Q_1(t,x)=d(\phi+V)-d(\phi)-d'(\phi)V \label{2.22-1}
\end{equation}
with $\phi=\phi(x_1+ct)$ and $V=V(t,x)$, and
\begin{equation}
Q_2(t-\tau,x-y)=b(\phi+V)-b(\phi)-b'(\phi)V \label{2.22}
\end{equation}
with $\phi=\phi(x_1-y_1+c(t-\tau))$ and $V=V(t-\tau,x-y)$. Here
$I_i, \ i=1,2,3,4$, denotes the $i$-th term in the right-side of
line above (\ref{2.21}).

From (H$_3$), i.e., $d''(u)\ge 0$ and $b''(u)\le 0$, applying
Taylor formula to \eqref{2.22-1} and \eqref{2.22}, we immediately
have
\[
Q_1(t,x)\ge 0 \ \mbox{ and } \ Q_2(t-\tau,x-y)\le 0,
\]
which implies
\begin{equation}
I_1(t,x)\le 0 \ \mbox{ and } \ I_2(t-\tau,x)\le 0. \label{2.22-new0}
\end{equation}
From (H$_3$) again, since $d'(\phi)$ is increasing and $b'(\phi)$ is
decreasing, then  $d'(0)-d'(\phi(x_1+ct))\le 0$ and
$b'(\phi(x_1-y_1+c(t-\tau)))-b'(0)\le 0$, which imply, with $V\ge 0$,
\begin{equation}
I_3(t,x)\le 0 \ \mbox{ and } \ I_4(t-\tau,x)\le 0. \label{2.22-new1}
\end{equation}
Thus, applying \eqref{2.22-new0} and \eqref{2.22-new1} to
\eqref{2.21}, we obtain
\begin{equation}
\dfrac{\partial V}{\partial t}-J\ast V+V+d'(0) V-b'(0)
\int_{\mathbb{R}^n} f_\beta(y) V(t-\tau,x-y)dy\le 0. \label{2.23}
\end{equation}

Let ${\bar V}(t,x)$ be the solution of the following equation with
the same initial data  $V_0(s,x)$:
\begin{equation} \label{2010-10}
\begin{cases}
\dfrac{\partial {\bar V}}{\partial t}-J\ast {\bar V}+\bar{V} +d'(0)
{\bar V} -  b'(0) \displaystyle \int_{\mathbb{R}^n}
f_\beta(y) {\bar V}(t-\tau,x-y) dy = 0,\ \ \ \
 (t,x)\in R_+\times \mathbb{R}^n,\\
{\bar V}(s,x)=V_0(s,x), \ \ \ s\in[-\tau,0], x\in \mathbb{R}^n.
\end{cases}
\end{equation}
From Proposition \ref{lemma6}, we know that $\bar{V}(t,x)$ globally
exists. Furthermore, \eqref{2010-10} is actually a linear equation, and its solution
 is as smooth as its  initial data.
By the comparison principle (Proposition \ref{lemma7}), we have
\begin{equation}
0\leq V(t,x) \le {\bar V}(t,x), \ \ \ \mbox{ for } (t,x)\in
\mathbb{R}_+\times \mathbb{R}^n. \label{2010-11}
\end{equation}

Let
\begin{equation}
v(t,x):=e^{-\lambda_*(x_1+ct-x_*)}{\bar V}(t,x). \label{2010-12}
\end{equation}
From \eqref{2010-10}, $v(t,x)$ satisfies
\begin{eqnarray}
&&\dfrac{\partial {v}}{\partial
t}-\int_{\mathbb{R}^n}J(y)e^{-\lambda_\ast y_1}v(t,x-y)dy +c_1
v\nonumber\\
&&\qquad\qquad = c_2 \int_{\mathbb{R}^n} f_\beta(y)
e^{-\lambda_*(y_1+c\tau)}v(t-\tau,x-y) dy , \label{2010-13-new-new}
\end{eqnarray}
where
\begin{equation}
 c_1:=c\lambda_*+1+d'(0)>0, \
\mbox{ and } \ c_2:=b'(0). \label{2010-14}
\end{equation}
When $\tau=0$, then \eqref{2010-13-new} is reduced to
\begin{equation}
\dfrac{\partial {v}}{\partial t}
-\int_{\mathbb{R}^n}J(y)e^{-\lambda_\ast y_1}v(t,x-y)dy +c_1 v = c_2
\int_{\mathbb{R}^n} f_\beta(y) e^{-\lambda_* y_1}v(t,x-y) dy .
\label{2010-13-2}
\end{equation}
Applying Proposition \ref{lemma4} to \eqref{2010-13-new-new} for $\tau>0$ and Proposition \ref{lemma5} to \eqref{2010-13-2} for $\tau=0$, we obtain the following
decay rates:
\begin{eqnarray}
& &\|v(t)\|_{L^\infty(\mathbb{R}^n)}\le C t^{-\frac{n}{\alpha}}
e^{-\varepsilon_1(\tilde{c}_1-c_3)t},   \ \ \ \mbox{ for } \ \tau>0, \label{new-1} \\
& &\|v(t)\|_{L^\infty(\mathbb{R}^n)}\le C t^{-\frac{n}{\alpha}}
e^{-(\tilde{c}_1-c_3)t},   \ \ \ \mbox{ for } \ \tau=0, \label{new-1-new}
\end{eqnarray}
where $0<\varepsilon_1=\varepsilon_1(\tau)<1$, and  $c_3$ is defined in \eqref{p21}, which can be directly
calculated as, by using the property \eqref{2.2-3},
\begin{eqnarray}
c_3&=&b'(0)\int_{\mathbb{R}^n} f_\beta(y) e^{-\lambda_*(y_1+c\tau)} dy \notag \\
&=&b'(0)\int_{\mathbb{R}} f_{1\beta}(y_1) e^{-\lambda_*(y_1+c\tau)} dy_1 \notag \\
&=&b'(0)e^{\beta\lambda_*^2-\lambda_* c \tau}>0. \label{new-2}
\end{eqnarray}
and
\begin{eqnarray}
\tilde{c}_1=c\lambda_*+1+d'(0)-\int_{\mathbb{R}}J(y_1)e^{-\lambda_\ast
y_1}dy_1=c\lambda_*+d'(0)-E_c(\lambda_\ast). \label{new-2-2}
\end{eqnarray}

When $c>c_*$, namely, the wave $\phi(x_1+ct)$ is non-critical, from
\eqref{2.4'-newnew} in Theorem \ref{TW}, we realize
\begin{equation}
\tilde{c}_1:=c\lambda_*+d'(0)-E_c(\lambda_\ast)=G_c(\lambda_*)>H_c(\lambda_*)=b'(0)e^{\beta
\lambda_*^2-\lambda_* c\tau}=:c_3. \label{new-3}
\end{equation}
Thus, \eqref{new-1} and \eqref{new-1-new}  immediately imply the following exponential
decay for $c>c_*$
\begin{eqnarray}
& &\|v(t)\|_{L^\infty(\mathbb{R}^n)}\le C t^{-\frac{n}{\alpha}}
e^{-\varepsilon_1\tilde{\mu} t}, \ \ \mbox{ for } \tau>0, \label{new-4} \\
& &\|v(t)\|_{L^\infty(\mathbb{R}^n)}\le C t^{-\frac{n}{\alpha}}
e^{-\tilde{\mu} t}, \ \ \mbox{ for } \tau=0, \label{new-4-new}
\end{eqnarray}
where
\begin{equation}
\tilde{\mu}:=\tilde{c}_1-c_3=G_c(\lambda_*)-H_c(\lambda_*)>0. \label{mu-1}
\end{equation}
When $c=c_*$, namely, the wave $\phi(x_1+c_*t)$ is critical, from
\eqref{2.5-new} in Proposition \ref{TW}, we realize
\begin{equation}
\tilde{c}_1:=c\lambda_*+d'(0)-E_c(\lambda_\ast)=G_c(\lambda_*)=H_c(\lambda_*)=b'(0)e^{\beta
\lambda_*^2-\lambda_* c\tau}:=c_3. \label{new-5}
\end{equation}
Then, from \eqref{new-1} and \eqref{new-1-new}, we immediately obtain the following
algebraic decay for $c=c_*$
\begin{equation}
\|v(t)\|_{L^\infty(\mathbb{R}^n)}\le C t^{-\frac{n}{\alpha}}, \ \  \ \mbox{ for all} \tau\ge 0.
\label{new-6}\\
\end{equation}

Since $V(t,x)\leq {\bar V}(t,x)=e^{\lambda_*(x_1+ct-x_*)}v(t,\xi)$,
and $0<e^{\lambda_*(x_1+ct-x_*)}\leq 1$ for $x_1\in
(-\infty,x_{*}-ct]$, we immediately obtain the following decay for
$V$.

\begin{lemma} \label{new-lemma3}
Let $V=V(t,x)$. Then
\begin{enumerate}
\item when $c>c_*$, then
\begin{eqnarray}
& & \|V(t)\|_{L^\infty((-\infty, \ x_{*}-ct]\times\mathbb{R}^{n-1})} \leq C(1+t)^{-\frac{n}{\alpha}} e^{-\varepsilon_1\tilde{\mu} t},
\ \ \ \mbox{ for } \tau>0,\\
& & \|V(t)\|_{L^\infty((-\infty, \ x_{*}-ct]\times\mathbb{R}^{n-1})} \leq C(1+t)^{-\frac{n}{\alpha}} e^{-\tilde{\mu} t},
\ \ \ \mbox{ for } \tau=0;
\label{2010-27}
\end{eqnarray}
Here $\tilde{\mu}:=\tilde{c}_1-c_3=G_c(\lambda_*)-H_c(\lambda_*)>0$ for $c>c_*$.
\item when $c=c_*$, then
\begin{equation}
 \|V(t)\|_{L^\infty((-\infty, \ x_{*}-ct]\times\mathbb{R}^{n-1})} \le C(1+t)^{-\frac{n}{\alpha}}, \ \ \ \mbox{ for all } \tau\ge 0.
\label{2010-28}
\end{equation}
\end{enumerate}
\end{lemma}

Next we prove $V(t,x)$ exponentially decay for $x \in
[x_*-ct,\infty) \times \mathbb{R}^{n-1}$.

\begin{lemma}\label{new-lemma4} For   $\tau>0$,  it holds that
\begin{eqnarray}
&\|V(t)\|_{L^\infty([x_*-ct,\infty)\times \mathbb{R}^{n-1})} \leq C
t^{-\frac{n}{\alpha}}e^{-\mu_{\tau} t}, & \ \mbox{for}\  c>c_*
\label{new-7},\\
&\|V(t)\|_{L^\infty([x_*-ct,\infty)\times \mathbb{R}^{n-1})} \leq C
t^{-\frac{n}{\alpha}}, & \ \mbox{for}\  c=c_* \label{new-7-2},
\end{eqnarray}
with some constant
$0<\mu_{\tau}<\min\{d'(u_+)-b'(u_+),\varepsilon_1\tilde{\mu}\}$ for $c>c_*$.
\end{lemma}

\noindent {\bf Proof}. From \eqref{2.16-1} and \eqref{2.2}, as set
in \eqref{2.19} $V(t,x):=U^+(t,x)-\phi(x_1+ct)$, we have
\begin{equation}
\dfrac{\partial V}{\partial t} -J\ast V+V +d(\phi+V)-d(\phi) =
\int_{\mathbb{R}^n} f_\beta(y) [b(\phi+V)-b(\phi)] dy.
 \label{new-7-1}
\end{equation}
Applying Taylor expansion formula and noting (H$_3$) for $d''(u)\ge
0$ and $b''(u)\le 0$, we have
\begin{eqnarray}
& &d(\phi+V)-d(\phi)=d'(\phi)V+d''({\bar\phi}_1)V^2\ge d'(\phi)V, \label{nov1}\\
& &b(\phi+V)-b(\phi)=b'(\phi)V+b''({\bar\phi}_2)V^2\le b'(\phi)V,
\label{nov2}
\end{eqnarray}
where ${\bar\phi}_i$ ($i=1,2$)  are some functions between $\phi$
and $\phi+V$. Substituting \eqref{nov1} and \eqref{nov2} into
\eqref{new-7-1}, and noticing Lemma \ref{new-lemma3},  we have
\begin{equation}\label{nov3}
\begin{cases}
\displaystyle \dfrac{\partial V}{\partial t}
-J\ast V+V +d'(\phi)V  \le  \int_{\mathbb{R}^n} f_\beta(y) b'(\phi(x_1-y_1+c(t-\tau))) V(t-\tau,x-y) dy,  \\
\qquad \qquad \qquad \qquad \qquad \qquad \qquad \qquad \quad \
\mbox{ for } t>0, x\in  \mathbb{R}^{n}\\
V|_{x_1\leq x_*-ct}\le C_2 (1+t)^{-\frac{n}{\alpha}}e^{-\varepsilon_1\tilde{\mu} t},  \qquad  \ \ \mbox{ for } t>0, (x_2,\cdots,x_n)\in  \mathbb{R}^{n-1}\\
V|_{t=s}=V_0(s,x),  \qquad \qquad \qquad \qquad \quad \ \ \ \mbox{
for } s\in [-\tau,0], x\in \mathbb{R}^{n}
\end{cases}
\end{equation}
for some positive constant $C_2$.

Let
\begin{equation}
\tilde{V}(t)=C_3 (1+\tau+t)^{-\frac{n}{\alpha}}e^{-\mu_{\tau} t}
\label{nov5}
\end{equation}
for $C_3\ge C_2\ge \max_{(s,x)\in[-\tau,0]\times\mathbb{R}^n}|V_0(s,x)|$.
As  in \eqref{nov4-new}, for given $0<\varepsilon_0 <1$, we can select a sufficiently large number  $x_*$  such that, for
$\xi_1\ge x_*\gg 1$,
\begin{equation}
d'(\phi(\xi_1))-\int_{\mathbb{R}^n} f_\beta(y)
b'(\phi(\xi_1-y_1-c\tau)) dy\geq \varepsilon_0 [d'(u_+)-b'(u_+)]>0.
\label{nov4}
\end{equation}
Thus, we have
\begin{eqnarray} \label{01-20}
& & \displaystyle \frac{\partial \tilde{V}}{\partial t} -J\ast\tilde{V} + \tilde{V}
+d'(\phi)\tilde{V} - \int_{\mathbb{R}^n} f_\beta(y)
b'(\phi(\xi_1-y_1-c\tau))
\tilde{V}(t-\tau ) dy \notag \\
& &=-\frac{n}{\alpha}C_3(1+t+\tau)^{-\frac{n}{\alpha}-1} e^{-\mu_{\tau} t} -\mu_{\tau} C_3(1+t+\tau)^{-\frac{n}{\alpha}} e^{-\mu_{\tau} t} \notag \\
& &\ \ \ +C_3(1+t+\tau)^{-\frac{n}{\alpha}} e^{-\mu_{\tau} t}d'(\phi(\xi_1)) \notag \\
& &\ \ \ -C_3(1+t)^{-\frac{n}{\alpha}} e^{-\mu_{\tau} (t-\tau)} \int_{\mathbb{R}^n} f_\beta (y) b'(\phi(\xi_1-y_1-c\tau))dy \notag \\
& &=C_3(1+t+\tau)^{-\frac{n}{\alpha}} e^{-\mu_{\tau} t}\Big\{\Big[d'(\phi(\xi_1))-\int_{\mathbb{R}^n} f_\beta (y) b'(\phi(\xi_1-y_1-c\tau))dy \Big]-\mu_{\tau}
\notag \\
& & \ \ \  -\frac{n}{\alpha}(1+t+\tau)^{-1} -\left(e^{\mu_{\tau} \tau}\left(\frac{1+t}{1+t+\tau}\right)^{-\frac{n}{\alpha}}-1\right) \int_{\mathbb{R}^n} f_\beta (y) b'(\phi(\xi_1-y_1-c\tau))dy\Big\}
\notag \\
& &\ge C_3(1+t+\tau)^{-\frac{n}{\alpha}} e^{-\mu_{\tau} t}\Big\{\varepsilon_0 [d'(u_+)-b'(u_+)] -\mu_{\tau} -\frac{n}{\alpha}(1+t+\tau)^{-1}
\notag \\
& & \ \ \ -\left(e^{\mu_{\tau} \tau}\left(\frac{1+t}{1+t+\tau}\right)^{-\frac{n}{\alpha}}-1\right) \int_{\mathbb{R}^n} f_\beta (y) b'(\phi(\xi_1-y_1-c\tau))dy\Big\} \notag \\
& &\ge 0
\label{01-21}
\end{eqnarray}
by selecting a sufficiently small number
\begin{eqnarray}
&0<\mu_{\tau} <d'(u_+)-b'(u_+) &  \mbox{ for } c>c_*, \label{01-23}\\
&\mu_{\tau}=0  & \mbox{ for } c=c_*, \label{01-24}
\end{eqnarray}
 and taking $t\ge l_0\tau$ for a sufficiently large integer $l_0\gg 1$. Hence, we proved that
\begin{eqnarray}\label{nov6}
\begin{cases}
\displaystyle \frac{\partial \tilde{V}}{\partial t} -J\ast\tilde{V} + \tilde{V}
+d'(\phi)\tilde{V} \geq  \int_{\mathbb{R}^n} f_\beta(y)
b'(\phi(\xi_1-y_1-c\tau))
\tilde{V}(t-\tau) dy,\\
\qquad\qquad\qquad\qquad\qquad\qquad\qquad\qquad\qquad \mbox{ for }
t>l_0\tau, \xi\in [x_*,+\infty)\times \mathbb{R}^{n-1}\\
\tilde{V}|_{\xi_1= x_*}=C_3
(1+\tau+t)^{-\frac{n}{\alpha}}e^{-\mu_{\tau} t} > C_2
(1+t)^{-\frac{n}{\alpha}}e^{-\varepsilon_1\tilde{\mu} t},
\ \mbox{ for } t>0, (\xi_2,\cdots,\xi_n)\in  \mathbb{R}^{n-1}\\
\tilde{V}(t)=C_3 (1+\tau+t)^{-\frac{n}{\alpha}}e^{-\mu_{\tau}
t} > V_0(t,\xi),  \qquad \ \ \ \mbox{ for } t\in
[-\tau,l_0\tau], \xi\in \mathbb{R}^{n}.
\end{cases}
\end{eqnarray}

Denote $\Omega:=\{(x,t)|x_1\geq x_*-ct,\ t\geq l_0\tau,\
(x_2,\cdots,x_n)\in  \mathbb{R}^{n-1}\}$. Noticing the construction
of \eqref{nov3} and \eqref{nov6},  then similar to the proof of
Proposition \ref{lemma7} , we know that
\begin{equation}\label{4.15}
\tilde{V}(t)- V(t,x)\geq 0,\ \mbox{for}\
(x,t)\in\mathbb{R}^n\times[-\tau,\infty)\setminus \Omega.
\end{equation}
Thus the proof is complete. $\square$

\medskip

For $\tau=0$, it is easy to prove the corresponding results as follows.

\begin{lemma}\label{new-lemma4-new} For   $\tau=0$,  it holds that
\begin{eqnarray}
&\|V(t)\|_{L^\infty([x_*-ct,\infty)\times \mathbb{R}^{n-1})} \leq C
t^{-\frac{n}{\alpha}}e^{-\mu_{\tau} t}, & \ \mbox{for}\  c>c_*
\label{new-7-new},\\
&\|V(t)\|_{L^\infty([x_*-ct,\infty)\times \mathbb{R}^{n-1})} \leq C
t^{-\frac{n}{\alpha}}, & \ \mbox{for}\  c=c_* \label{new-7-2-new},
\end{eqnarray}
with some constant
$0<\mu_{\tau}<\min\{d'(u_+)-b'(u_+),\varepsilon_1\tilde{\mu}\}$ for $c>c_*$.
\end{lemma}

Combing Lemma \ref{new-lemma3}-Lemma \ref{new-lemma4-new}, we obtain
the decay rates for $V(t,x)$ in $L^\infty(\mathbb{R}^n)$.

\begin{lemma} \label{new-lemma5} It holds that:
\begin{enumerate}
\item when $c>c_*$, then
\begin{eqnarray}
& & \|V(t)\|_{L^\infty(\mathbb{R}^{n})} \le C(1+t)^{-\frac{n}{\alpha}} e^{-\mu_\tau t}, \quad \mbox{ for } \tau>0,
\label{2010-27-1} \\
& & \|V(t)\|_{L^\infty(\mathbb{R}^{n})} \le C(1+t)^{-\frac{n}{\alpha}} e^{-\mu_0 t}, \quad \mbox{ for } \tau=0,
\label{2010-27-1-new}
\end{eqnarray}
where $0<\mu_\tau<\min\{d'(u_+)-b'(u_+),\varepsilon_1[G_c(\lambda_*)-H_c(\lambda_*)]\}$ with $0<\varepsilon_1<1$ for
$\tau>0$, and $0<\mu_0< \min\{d'(u_+)-b'(u_+),G_c(\lambda_*)-H_c(\lambda_*)\}$ for $\tau=0$;
\item when $c=c_*$,
\begin{equation}
 \|V(t)\|_{L^\infty(\mathbb{R}^{n})} \le C(1+t)^{-\frac{n}{\alpha}}, \quad \mbox{ for all } \tau\ge 0.
\label{2010-28-1}
\end{equation}
\end{enumerate}
\end{lemma}

Since $V(t,x)=U^+(t,x)-\phi(x_1+ct)$, Lemma \ref{new-lemma5}  give
directly the following convergence for the solution in the cases
with time-delay.

\begin{lemma}\label{new-lemma6}
It holds that:
\begin{enumerate}
\item when $c>c_*$, then
\begin{eqnarray}
& & \sup_{x\in \mathbb{R}^{n}}|U^+(t,x)-\phi(x_1+ct)| \le C(1+t)^{-\frac{n}{\alpha}} e^{-\mu_\tau t}, \quad \mbox{ for } \tau>0,
\label{2010-27-2} \\
& & \sup_{x\in \mathbb{R}^{n}}|U^+(t,x)-\phi(x_1+ct)| \le C(1+t)^{-\frac{n}{\alpha}} e^{-\mu_0 t}, \quad \mbox{ for } \tau=0,
\label{2010-27-2-new}
\end{eqnarray}
where $0<\mu_\tau<\min\{d'(u_+)-b'(u_+),\varepsilon_1[G_c(\lambda_*)-H_c(\lambda_*)]\}$ with $0<\varepsilon_1<1$ for
$\tau>0$, and $0<\mu_0< \min\{d'(u_+)-b'(u_+),G_c(\lambda_*)-H_c(\lambda_*)\}$ for $\tau=0$;
\item when $c=c_*$, then
\begin{equation}
 \sup_{x\in \mathbb{R}^{n}}|U^+(t,x)-\phi(x_1+c_*t)|  \le C(1+t)^{-\frac{n}{\alpha}}, \quad \mbox{ for all } \tau\ge 0.
\label{2010-28-2}
\end{equation}
\end{enumerate}
\end{lemma}

{\bf Step 2. The convergence of $U^-(t,x)$ to $\phi(x_1+ct)$}

For the traveling wave $\phi(x_1+ct)$ with  $c\ge c_*$, let
\begin{equation}
V(t,x)=\phi(x_1+ct)-U^-(t,x), \ \ \
V_0(s,x)=\phi(x_1+cs)-U^-_0(s,x). \label{2.50}
\end{equation}
As in Step 1, we can similarly prove that $U^-(t,x)$ converges to
$\phi(x_1+ct)$ as follows.

\begin{lemma}\label{new-lemma7}
It holds that:
\begin{enumerate}
\item when $c>c_*$, then
\begin{eqnarray}
& & \sup_{x\in \mathbb{R}^{n}}|U^-(t,x)-\phi(x_1+ct)| \le C(1+t)^{-\frac{n}{\alpha}} e^{-\mu_\tau t}, \quad \mbox{ for } \tau>0,
\label{new2010-27-2} \\
& & \sup_{x\in \mathbb{R}^{n}}|U^-(t,x)-\phi(x_1+ct)| \le C(1+t)^{-\frac{n}{\alpha}} e^{-\mu_0 t}, \quad \mbox{ for } \tau=0,
\label{new2010-27-2-new}
\end{eqnarray}
where $0<\mu_\tau<\min\{d'(u_+)-b'(u_+),\varepsilon_1[G_c(\lambda_*)-H_c(\lambda_*)]\}$ with $0<\varepsilon_1<1$ for
$\tau>0$, and $0<\mu_0< \min\{d'(u_+)-b'(u_+),G_c(\lambda_*)-H_c(\lambda_*)\}$ for $\tau=0$;
\item when $c=c_*$, then
\begin{equation}
 \sup_{x\in \mathbb{R}^{n}}|U^-(t,x)-\phi(x_1+c_*t)|  \le C(1+t)^{-\frac{n}{\alpha}}, \quad \mbox{ for all } \tau\ge 0.
\label{new2010-28-2}
\end{equation}
\end{enumerate}
\end{lemma}

{\bf Step 3. The convergence of $u(t,x)$ to $\phi(x_1+ct)$}

Finally, we prove that $u(t,x)$ converges to $\phi(x_1+ct)$. Since
the initial data satisfy $U^-_0(s,x)\le u_0(s,x) \le U^+_0(s,x)$ for
$(s,x)\in [-\tau,0]\times \mathbb{R}^n$,  then the comparison principle implies that
\[
U^-(t,x)\le u(t,x) \le U^+(t,x), \ \ \ \ (t,x)\in R_+\times
\mathbb{R}^n.
\]
Thanks to Lemmas \ref{new-lemma6} and \ref{new-lemma7}, by the
squeeze argument, we  have the following convergence results.

\begin{lemma}\label{new-lemma8}
It holds that:
\begin{enumerate}
\item when $c>c_*$, then
\begin{eqnarray}
& & \sup_{x\in \mathbb{R}^{n}}|u(t,x)-\phi(x_1+ct)| \le C(1+t)^{-\frac{n}{\alpha}} e^{-\mu_\tau t}, \quad \mbox{ for } \tau>0,
\label{nn2010-27-2} \\
& & \sup_{x\in \mathbb{R}^{n}}|u(t,x)-\phi(x_1+ct)| \le C(1+t)^{-\frac{n}{\alpha}} e^{-\mu_0 t}, \quad \mbox{ for } \tau=0,
\label{nn2010-27-2-new}
\end{eqnarray}
where $0<\mu_\tau<\min\{d'(u_+)-b'(u_+),\varepsilon_1[G_c(\lambda_*)-H_c(\lambda_*)]\}$ with $0<\varepsilon_1<1$ for
$\tau>0$, and $0<\mu_0< \min\{d'(u_+)-b'(u_+),G_c(\lambda_*)-H_c(\lambda_*)\}$ for $\tau=0$;
\item when $c=c_*$, then
\begin{equation}
 \sup_{x\in \mathbb{R}^{n}}|u(t,x)-\phi(x_1+c_*t)|  \le C(1+t)^{-\frac{n}{\alpha}}, \quad \mbox{ for all } \tau\ge 0.
\label{nn2010-28-2}
\end{equation}
\end{enumerate}
\end{lemma}

\section{Applications and Concluding Remark}

In this section, we first give the direct applications of Theorem \ref{TW}-\ref{thm1} to the Nicholson's blowflies type
equation with nonlocal dispersion,  and the classical Fisher-KPP equation with nonlocal dispersion. Then we point out that,
the developed stability theory above can be also applied to the more general case.

\subsection{Nicholson's blowflies  equation with nonlocal dispersion}

For the equation \eqref{1.1}, by taking $d(u)=\delta u$ and $b(u)=pue^{-au}$ with $\delta>0$, $p>0$ and $a>0$, we get the so-called Nicholson's blowflies equation with nonlocal dispersion

\begin{equation}
\begin{cases}
\displaystyle \frac{\partial u}{\partial t} - J\ast u +u+\delta u(t,x))=
p\int_{\mathbb{R}^n} f_\beta
(y) u(t-\tau,x-y)e^{-au(t-\tau,x-y)} dy, \\
u(s,x)=u_0(s,x), \ \ s\in[-\tau,0], \ x\in \mathbb{R}^n.
\end{cases}
\label{6.1}
\end{equation}
Clearly, there exist two constant equilibria $u_-=0$ and $u_+=\frac{1}{a}\ln \frac{p}{\delta}$, and the selected $d(u)$ and $b(u)$ satisfy the
hypothesis (H$_1$)-(H$_3$) automatically under the consideration of $1<\frac{p}{\delta}\le e$. Let $J(x)$ satisfy the hypothesis (J$_1$) and (J$_2$), from Theorem \ref{TW} and Theorem \ref{thm1}, we have
the following existence of monostable traveling waves and their stabilities.

\begin{theorem}[Traveling waves] Let $J(x)$ satisfy (J$_1$) and (J$_2$).
For \eqref{6.1}, there exists the minimal speed $c_*>0$,  such that
\begin{enumerate}
\item[$\bullet$] when $c\ge c_*$, the planar traveling waves $\phi(x\cdot{\bf e}_1 +ct)$  exist uniquely (up to a shift);
\item[$\bullet$] when $c< c_*$, the planar traveling waves $\phi(x\cdot{\bf e}_1 +ct)$  do not exist;
\end{enumerate}
Here $c_*>0$ and $\lambda_*>0$ are determined by
\[
H_{c_*}(\lambda_*)=G_{c_*}(\lambda_*) \ \mbox{ and } \  H'_{c_*}(\lambda_*)=G'_{c_*}(\lambda_*),
\]
where
\[
H_c(\lambda)=pe^{\beta \lambda^2-\lambda c\tau} \ \mbox{ and } \  G_c(\lambda)=c\lambda -\int_\mathbb{R}J_1(y_1)e^{-\lambda y_1}dy_1 +1+\delta.
\]
Particularly, when $c>c_*$, then
$
H_c(\lambda_*)<G_c(\lambda_*).
$
\end{theorem}

\begin{theorem}[Stability of traveling waves] Let $J(x)$ satisfy (J$_1$) and (J$_2$), and
the initial data be $u_0-\phi\in C([-\tau,0]; H^m_w(\mathbb{R}^n)\cap L^1_w(\mathbb{R}^n))$ and $\partial_s(u_0-\phi)\in L^1([-\tau,0]; H^{m+1}_w(\mathbb{R}^n)\cap L^1_w(\mathbb{R}^n))$ with $m>\frac{n}{2}$, and $u_-\le u_0\le u_+$ for $(s,x)\in [-\tau,0]\times \mathbb{R}^n$.
Then the solution of \eqref{6.1}  uniquely exists and
satisfies:
\begin{enumerate}
\item[$\bullet$] when $c>c_*$,  then
\begin{equation}
 \sup_{x\in \mathbb{R}^n}|u(t,x)-\phi(x_1+ct)|\le
C(1+t)^{-\frac{n}{\alpha}}e^{-\mu_\tau t}, \ \ t> 0, \label{2.12-new}
\end{equation}
for $0<\mu_\tau < \min\{d'(u_+)-b'(u_+), \  \varepsilon_1[G_c(\lambda_*)-H_c(\lambda_*)]\}$,
 and $\varepsilon_1=\varepsilon_1(\tau)$ such that $0<\varepsilon_1<1$ for $\tau>0$ and $\varepsilon_1=1$ for $\tau=0$

\item[$\bullet$] when $c=c_*$,  then
\begin{equation}
\sup_{x\in \mathbb{R}^n}|u(t,x)-\phi(x_1+c_*t)|\le
C(1+t)^{-\frac{n}{\alpha}}, \ \ t> 0. \label{2.12-2-new}
\end{equation}
\end{enumerate}
\end{theorem}

\subsection{Fisher-KPP equation with nonlocal dispersion}

For the equation \eqref{1.1}, let  $d(u)=u^2$, $b(u)=u$ and the delay  $\tau=0$, and take the limit of \eqref{1.1} as $\beta\to 0^+$, we get the classical Fisher-KPP equation with nonlocal dispersion without time-delay
\begin{equation}
\begin{cases}
\displaystyle \frac{\partial u}{\partial t} - J\ast u +u=u(1-u) \\
u(0,x)=u_0(x), \ \  \ x\in \mathbb{R}^n.
\end{cases}
\label{6.2}
\end{equation}
Then we have the existence of the monostable traveling waves and their stabilities from Theorem \ref{TW} and Theorem \ref{thm1}.

\begin{theorem}[Traveling waves] Let $J(x)$ satisfy (J$_1$) and (J$_2$).
For \eqref{6.2}, there exists the minimal speed $c_*>0$,  such that
\begin{enumerate}
\item[$\bullet$] when $c\ge c_*$, the planar traveling waves $\phi(x\cdot{\bf e}_1 +ct)$  exist uniquely (up to a shift);
\item[$\bullet$] when $c< c_*$, the planar traveling waves $\phi(x\cdot{\bf e}_1 +ct)$  do not exist;
\end{enumerate}
Here
$c_*:=\lambda_*^{-1} \int_{\mathbb{R}}J_1(y_1)e^{-\lambda_* y_1}dy_1$,
and $\lambda_*>0$ is determined by
$\int_{\mathbb{R}} (1+\lambda_* y_1) J_1(y_1) e^{-\lambda_* y_1} dy_1=0$.
When $c>c_*$, then
$
H_c(\lambda_*)<G_c(\lambda_*),
$
where
$
H_c(\lambda_*)=1$  and
$G_c(\lambda_*)=c\lambda_*-\int_{\mathbb{R}}J_1(y_1)e^{-\lambda_* y_1} dy_1 + 1$.
\end{theorem}

\begin{theorem}[Stability of traveling waves] Let $J(x)$ satisfy (J$_1$) and (J$_2$), and
the initial data be $u_0-\phi\in C([-\tau,0]; H^m_w(\mathbb{R}^n)\cap L^1_w(\mathbb{R}^n))$ with $m>\frac{n}{2}$, and $u_-\le u_0\le u_+$ for $x \in \mathbb{R}^n$.
Then the solution of \eqref{6.2}  uniquely exists and
satisfies:
\begin{enumerate}
\item[$\bullet$] when $c>c_*$,  then
\begin{equation}
 \sup_{x\in \mathbb{R}^n}|u(t,x)-\phi(x_1+ct)|\le
C(1+t)^{-\frac{n}{\alpha}}e^{-\mu_0 t}, \ \ t> 0, \label{2.12-nnew}
\end{equation}
for $0<\mu_0 <  \min\{d'(u_+)-b'(u_+), \  G_c(\lambda_*)-H_c(\lambda_*)\}$;

\item[$\bullet$] when $c=c_*$,  then
\begin{equation}
\sup_{x\in \mathbb{R}^n}|u(t,x)-\phi(x_1+c_*t)|\le
C(1+t)^{-\frac{n}{\alpha}}, \ \ t> 0. \label{2.12-2-nnew}
\end{equation}
\end{enumerate}
\end{theorem}

\subsection{Concluding Remark}
Here we give a remark on the wave stability to the generalized equations with nonlocal dispersion.
Let us consider  a more general monostable equation with nonlocal dispersion
\begin{equation}
\begin{cases}
\displaystyle \frac{\partial u}{\partial t} - J\ast u +u+d(u(t,x))=
F\Big(\int_{\mathbb{R}^n} \kappa
(y) b(u(t-\tau,x-y)) dy \Big), \\
u(s,x)=u_0(s,x), \ \ s\in[-\tau,0], \ x\in \mathbb{R}^n,
\end{cases}
\label{6.3}
\end{equation}
where $J(x)$ satisfies (J$_1$) and (J$_2$) as mentioned before, and $F(\cdot)$, $d(u)$, $b(u)$ and $g(x)$ satisfy
\begin{enumerate}
\item[$(\mathcal{H}_1)$] There exist $u_-=0$ and $u_+>0$ such that $d(0)=b(0)=F(0)=0$, $d(u_+)=F(b(u_+))$, $d\in C^2[0,u_+]$, $b\in C^2[0,u_+]$ and $F\in C^2[0,b(u_+)]$;
\item[$(\mathcal{H}_2)$] $F'(0)b'(0)>d'(0)\ge 0$ and $0<F'(b(u_+))b'(u_+)<d'(u_+)$;
\item[$(\mathcal{H}_3)$] $d'(u)\ge 0$, $b'(u)\ge 0$, $d''(u)\ge 0$ and $b''(u)\le 0$ for $u\in [0,u_+]$;
\item[$(\mathcal{H}_4)$] $F'(u)\ge 0$ and $F''(u)\le 0$ for $u\in [0,b(u_+)]$;
\item[$(\mathcal{H}_5)$] $\kappa(x)$ is a smooth, positive and  radial kernel with $\int_{\mathbb{R}^n}\kappa(x) dx=1$ and
$\int_{\mathbb{R}^n}\kappa(x)e^{-\lambda x_1}  dx <+\infty$ for all $\lambda>0$.
\end{enumerate}
Then, by a similar calculation, we can prove the existence of the traveling waves $\phi(x_1+ct)$ for $c\ge c_*$, where $c_*>0$
is a specified minimal wave speed, and that the noncritical traveling waves with $c>c_*$ are exponentially stable and the critical waves
with $c=c_*$ are algebraically stable.

\begin{theorem}[Traveling waves] Assume that $(\mbox{J}_1)$-$(\mbox{J}_2)$ and $(\mathcal{H}_1)$-$(\mathcal{H}_5)$ hold.
For \eqref{6.3}, there exists the minimal speed $c_*>0$,  such that
\begin{enumerate}
\item[$\bullet$] when $c\ge c_*$, the planar traveling waves $\phi(x\cdot{\bf e}_1 +ct)$  exist uniquely (up to a shift);
\item[$\bullet$] when $c< c_*$, the planar traveling waves $\phi(x\cdot{\bf e}_1 +ct)$  do not exist;
\end{enumerate}
Here $c_*>0$ and $\lambda_*=\lambda_*(c_*)>0$ are determined by
\[
\mathcal{H}_{c_*}(\lambda_*)=\mathcal{G}_{c_*}(\lambda_*) \ \mbox{ and } \ \mathcal{H}'_{c_*}(\lambda_*)=\mathcal{G}'_{c_*}(\lambda_*),
\]
where
\[
\mathcal{H}_c(\lambda)=F'(0)b'(0)\int_{\mathbb{R}^n} e^{-\lambda y_1} \kappa(y) dy, \ \ \
\mathcal{G}_c(\lambda)=c\lambda -\int_{\mathbb{R}} J_1(y_1) e^{-\lambda y_1} dy_1 +1+d'(0).
\]
When $c>c_*$, then
\[
H_c(\lambda_*)<G_c(\lambda_*).
\]
\end{theorem}

\begin{theorem}[Stability of traveling waves] Assume that $(\mbox{J}_1)$-$(\mbox{J}_2)$ and $(\mathcal{H}_1)$-$(\mathcal{H}_5)$ hold.
 Let
the initial data be $u_0-\phi\in C([-\tau,0]; H^{m+1}_w(\mathbb{R}^n)\cap L^1_w(\mathbb{R}^n))$   and $\partial_s(u_0-\phi)\in L^1([-\tau,0]; H^{m+1}_w(\mathbb{R}^n)\cap L^1_w(\mathbb{R}^n))$ with $m>\frac{n}{2}$, and $u_-\le u_0\le u_+$ for $x \in \mathbb{R}^n$.
Then the solution of \eqref{6.3}  uniquely exists and
satisfies:
\begin{enumerate}
\item[$\bullet$] when $c>c_*$,  then
\begin{equation}
 \sup_{x\in \mathbb{R}^n}|u(t,x)-\phi(x_1+ct)|\le
C(1+t)^{-\frac{n}{\alpha}}e^{-\mu_0 t}, \ \ t> 0, \label{2.12-nnnew}
\end{equation}
for $0<\mu_\tau < \min\{d'(u_+)-b'(u_+), \ \varepsilon_1[ G_c(\lambda_*)-H_c(\lambda_*)]\}$, and $0<\varepsilon_1<1$
for $\tau>0$ and $\varepsilon_1=1$ for $\tau=0$;

\item[$\bullet$] when $c=c_*$,  then
\begin{equation}
\sup_{x\in \mathbb{R}^n}|u(t,x)-\phi(x_1+c_*t)|\le
C(1+t)^{-\frac{n}{\alpha}}, \ \ t> 0. \label{2.12-2-nnnew}
\end{equation}
\end{enumerate}
\end{theorem}

\noindent {\bf Acknowledgments} The research of MM was supported in part by
Natural Sciences and Engineering Research Council of Canada under
the NSERC grant RGPIN 354724-08. The research of RH was supported in part by
NNSFC (No. 11001103) and SRFDP (No. 200801831002).

\end{document}